\documentclass[leqno,11pt,a4paper]{article}
\usepackage[left=3.2cm, top=2.2cm,bottom=2.5cm]{geometry}
\usepackage{amsmath,amssymb,amsthm,mathrsfs,dsfont,calc,dsfont,graphicx}
\usepackage[british]{babel}

\newtheorem{theorem}{Theorem}

\newtheorem{remark}{Remark}

\newcommand{\la}{\lambda}

\newcommand{\NN}{{\mathbb{N}}}
\newcommand{\RR}{{\mathbb{R}}}
\newcommand{\EE}{{\mathbb{E}}}

\def\EE{\mathbb{E}}
\def\NN{\mathbb{N}}

\def\RR{\mathbb{R}}

\def\eps{\varepsilon}
\def\la{\lambda}

\def\L{\Lambda}

\def\bP{\mathbb{P}}
\def\bQ{\mathbb{Q}}

\def\cB{\mathcal{B}}
\def\cH{\mathcal{H}}

\def\cP{\mathcal{P}}

\def\sE{\mathbf{E}}
\def\sY{\mathbf{G}}
\def\sI{\mathbf{I}}
\def\sY{\mathbf{Y}}
\def\sP{\mathbf{P}}
\def\sX{\mathbf{X}}
\def\sX{\mathbf{X}}
\def\sY{\mathbf{Y}}
\def\sV{\mathbf{V}}
\def\sZ{\mathbf{Z}}

\begin{document}

\title{The combinatorial structure of spatial STIT tessellations}
\author{Christoph Th\"ale and Viola Wei\ss}
\date{}
\maketitle

\begin{abstract}
Spatially homogeneous random tessellations that are {\bf st}able under {\bf it}eration (nesting) in the $3$-dimensional Euclidean space are considered, so-called STIT tessellations. They arise as outcome of a spatio-temporal process of subsequent cell division and consequently they are not facet-to-facet. The intent of this paper is to develop a detailed analysis of the combinatorial structure of such tessellations and to determine a number of new geometric mean values, for example for the neighborhood of the typical vertex. The heart of the results is a fine classification of tessellation edges based on the type of their endpoints or on the equality relationship with other types of line segments. In the background of the proofs are delicate distributional properties of spatial STIT tessellations.  
\end{abstract}
\begin{flushleft}\footnotesize
\end{flushleft}

\section{Introduction}\label{sec:intro}

A random tessellation (or mosaic) of a $d$-dimensional Euclidean space is a locally finite family of pairwise non-overlapping $d$-dimensional random convex polytopes -- called cells of the tessellation -- that cover the whole space. Whereas non-random tessellations and tilings are a central object in discrete geometry, random tessellations form one of the classical random structures considered in stochastic geometry (see \cite{SW,SKM} for some classical theory and \cite{Calka,Heinrich,HugSchneider,Schneider} for some more recent developments). The most popular models are hyperplane or Voronoi tessellations, where mainly the Poisson case has been studied. All these tessellations share the property of being side-to-side or facet-to-facet in higher dimensions. In recent years there has been a growing interest also in tessellation models that do not fulfill this property. In \cite{WC} a first systematic study of the complications is given when a tessellation is not side-to-side or facet-to-facet. Tessellations of this kind arise for example by subsequent cell division, see \cite{Cowan10}, which makes them interesting for particular applications in geology, material sciences or biology. However, such models have a much more complex structure, because the spatial arrangement of their cells can be rather complicated. This implies that tessellations of this kind carry a rich combinatorial structure. For example it is no more necessarily the case that the number of facets of a cell coincides with the number of its neighboring cells. It is therefore interesting to calculate a number of geometric mean values for the cells, its facets, its ridges, its vertices and other interesting classes of objects that are determined by the tessellation. They describe in some sense the `mean' cell-arrangement and allow some insight to the complex tessellation geometry.

Among these models for subsequent cell division the iteration stable or STIT tessellations are of particular interest in recent time in stochastic geometry, because of the number of analytic available results, see \cite{Cowan11,LR,NW06,NW08,ST,TW,TWN} and the references cited therein. A detailed combinatorial analysis of planar STIT tessellations has been carried out in \cite{Cowan11,NW06} and the intent of this paper is to study the combinatorial structure of STIT tessellations in $\RR^3$. The structure of spatial STIT tessellations is very rich and considerably more complex compared with the planar case, as we will see below. This is for the STIT tessellations mainly due to the fact that they do not only have T-shaped vertices as in the planar case, but also vertices of so-called X-type. This causes that in a detailed analysis we have to take into account the effect generated by the different types of vertices. And, moreover, many distributional results depend on the joint distribution of T- and X-vertices in a delicate way. It is worth noticing that some basic mean values for the spatial STIT model have already been determined in \cite{NW08,TW}. They will be in the background of our analysis. Based on the different vertex types a rather fine classification of the edges of the tessellation is introduced, similar to the edge classification of planar STIT tessellations in \cite{Cowan11}. The main result of the present paper is the investigation of the proportion of these different edge types in the tessellation. These breakdowns of the class of edges into different subclasses are based on fine distributional properties for spatial STIT tessellation that were developed in \cite{TWN} and which will be further developed below. They finally allow us to calculate a couple of new mean values, which describe the fine combinatorial structure and considerably refine the results from \cite{NW08,TW}.

The paper is organized as follows. We recall the basic construction of STIT tessellations in Subsection \ref{subsec:STITs} and recall in Subsection \ref{subse:typicalobjects} the definition and interpretation of object intensities, typical objects and object adjacencies. In Section \ref{sec:basicstructure} we introduce the basic geometric objects determined by a spatial STIT tessellation and their basic combinatorial structure from \cite{NW08,TW} in order to keep the paper self-contained. Our main results are the content of Section \ref{sec:finestructure}. Before presenting them we first take care of some vertex geometry in Subsection \ref{subsec:vertex geometry} and develop as a technical tool some new identities for so-called I-segments in Subsection \ref{subsec:isegmentdistributions1} that play an important r\^ole in our further analysis. In Subsection \ref{subsec:edgeclassification} we explain our refined edge classifications and calculate the probabilities that the typical edge of a STIT tessellation belongs to one of our classes. These probabilities are needed in Subsection \ref{subsec:newmeanvalues} to express a number of new geometric mean values that are related to a fine combinatorial structure of spatial STIT tessellations. A summary and a discussion of our results are presented in Section \ref{sec:summary}, whereas all proofs are the content of the final Section \ref{sec:proofs}.

\section{Preliminaries}\label{sec:preliminaries}

This section is devoted to some preparatory work. We first recall the basic construction of spatial STIT tessellations and thereafter we formalize the concept of intensities as well as typical and adjacent objects of classes of geometric objects that are related to our model. Moreover, some basic notation will be fixed.

\subsection{STIT tessellations}\label{subsec:STITs}
\begin{figure}
\begin{center}
 \includegraphics[width=0.32\columnwidth]{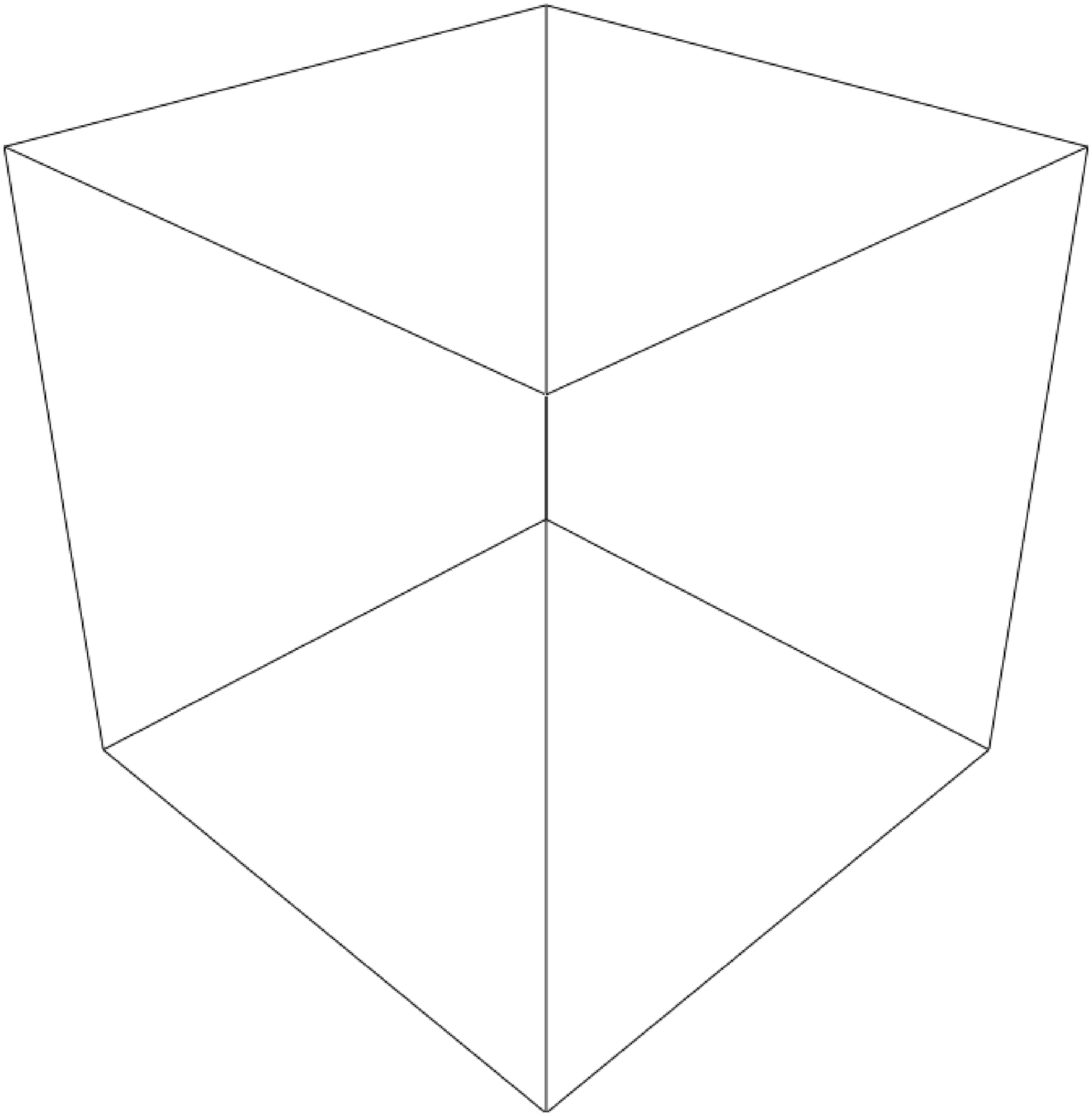}
 \includegraphics[width=0.32\columnwidth]{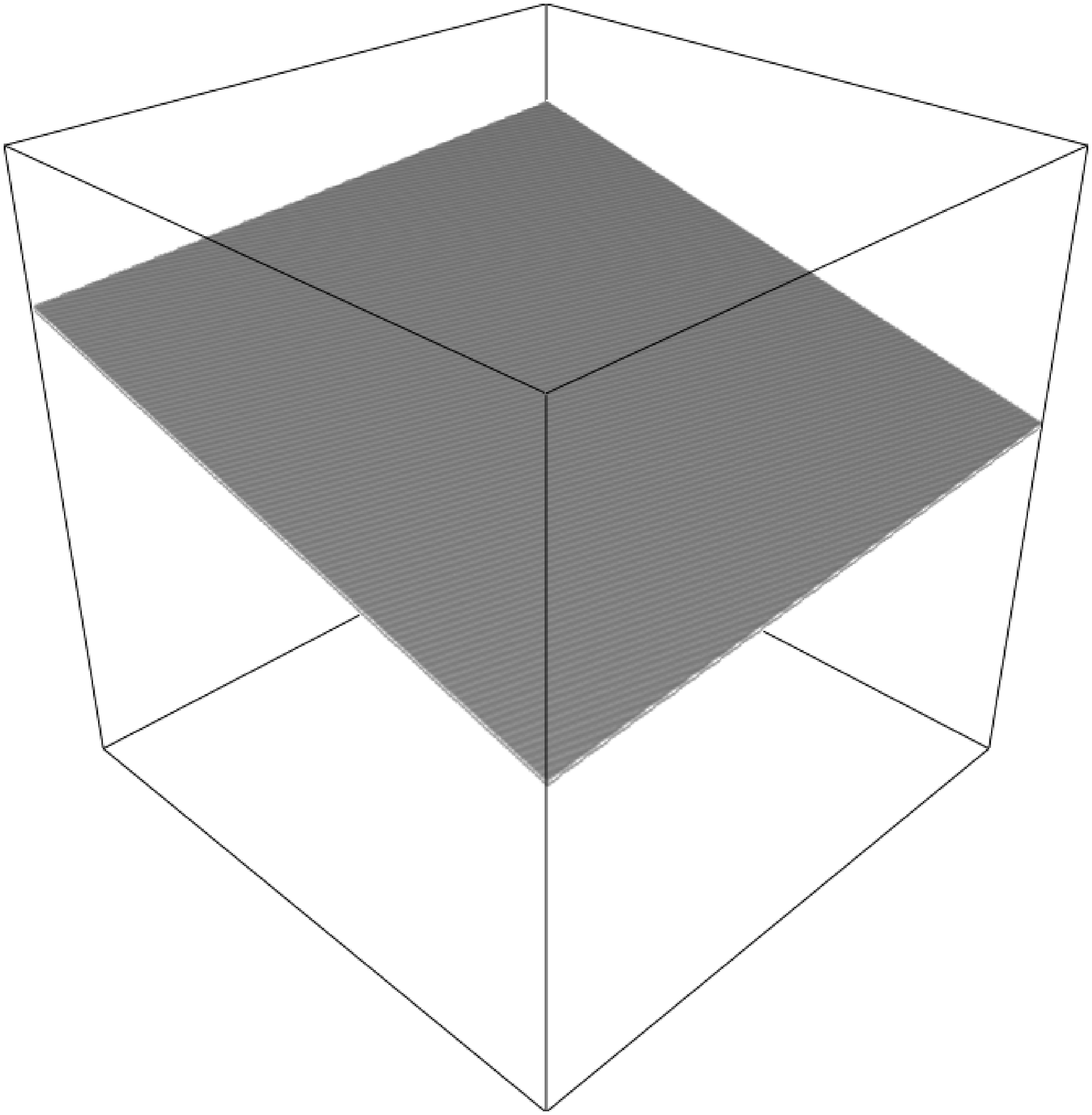}
 \includegraphics[width=0.32\columnwidth]{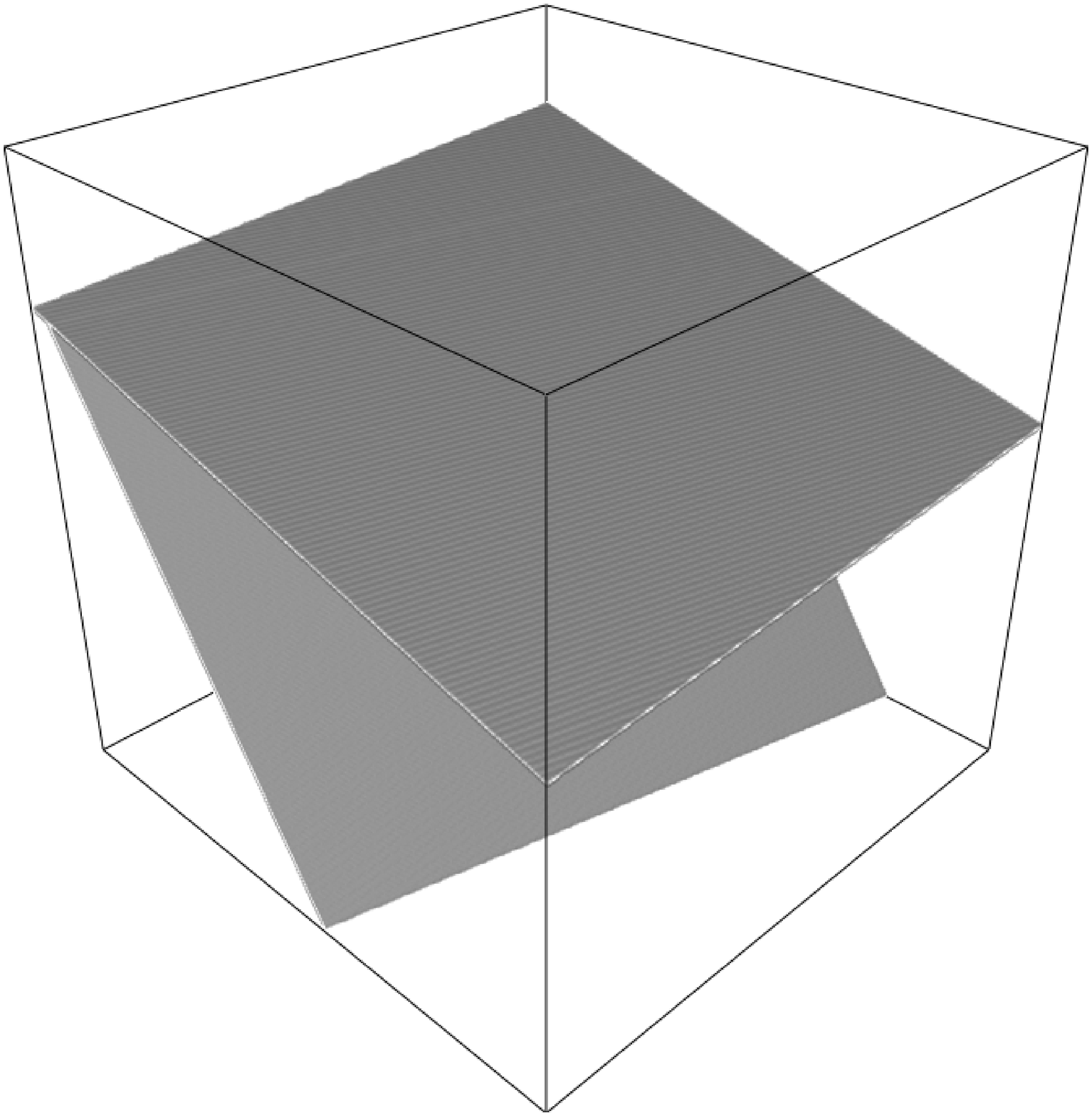}\\
 \includegraphics[width=0.32\columnwidth]{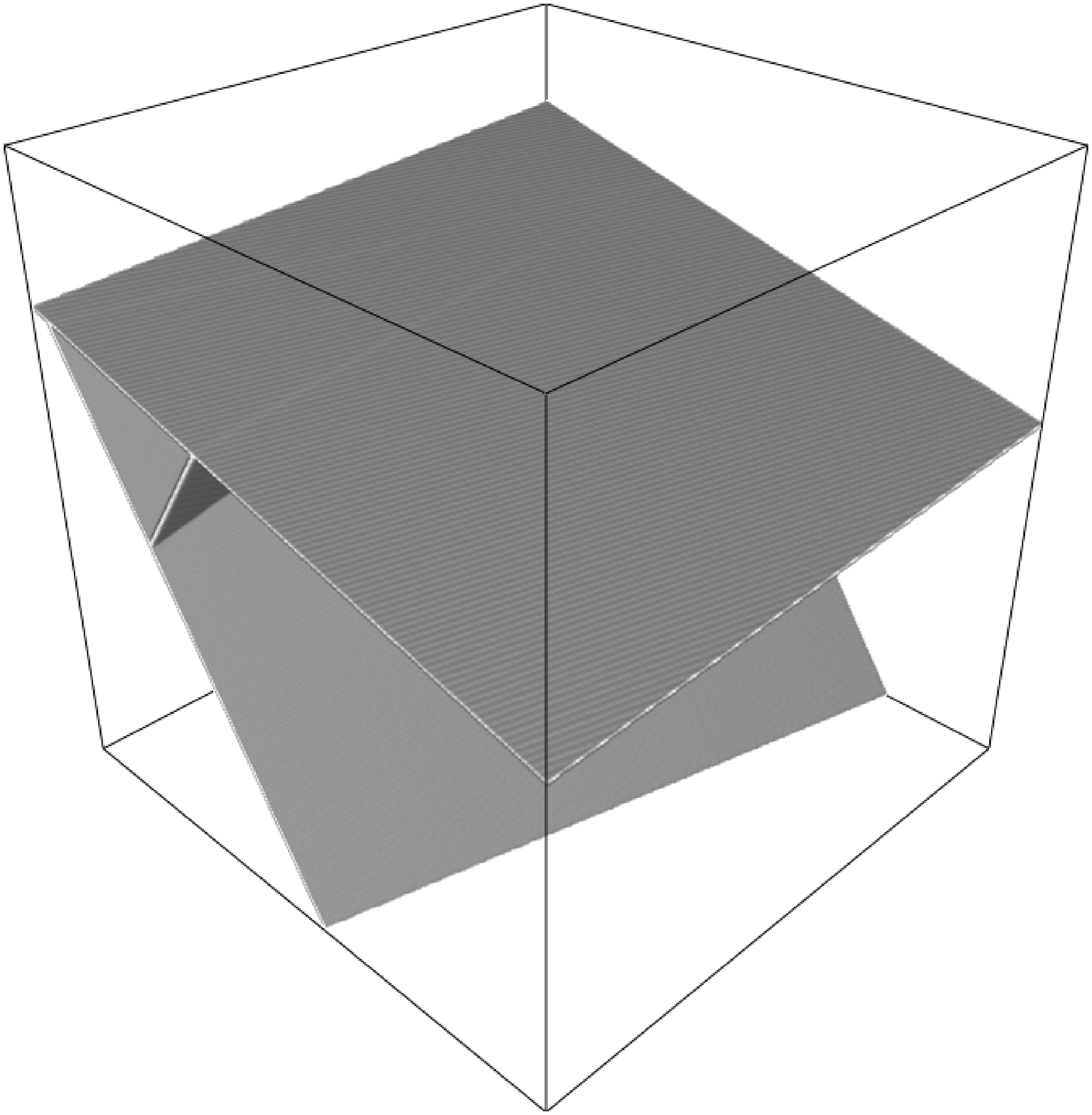}
 \includegraphics[width=0.32\columnwidth]{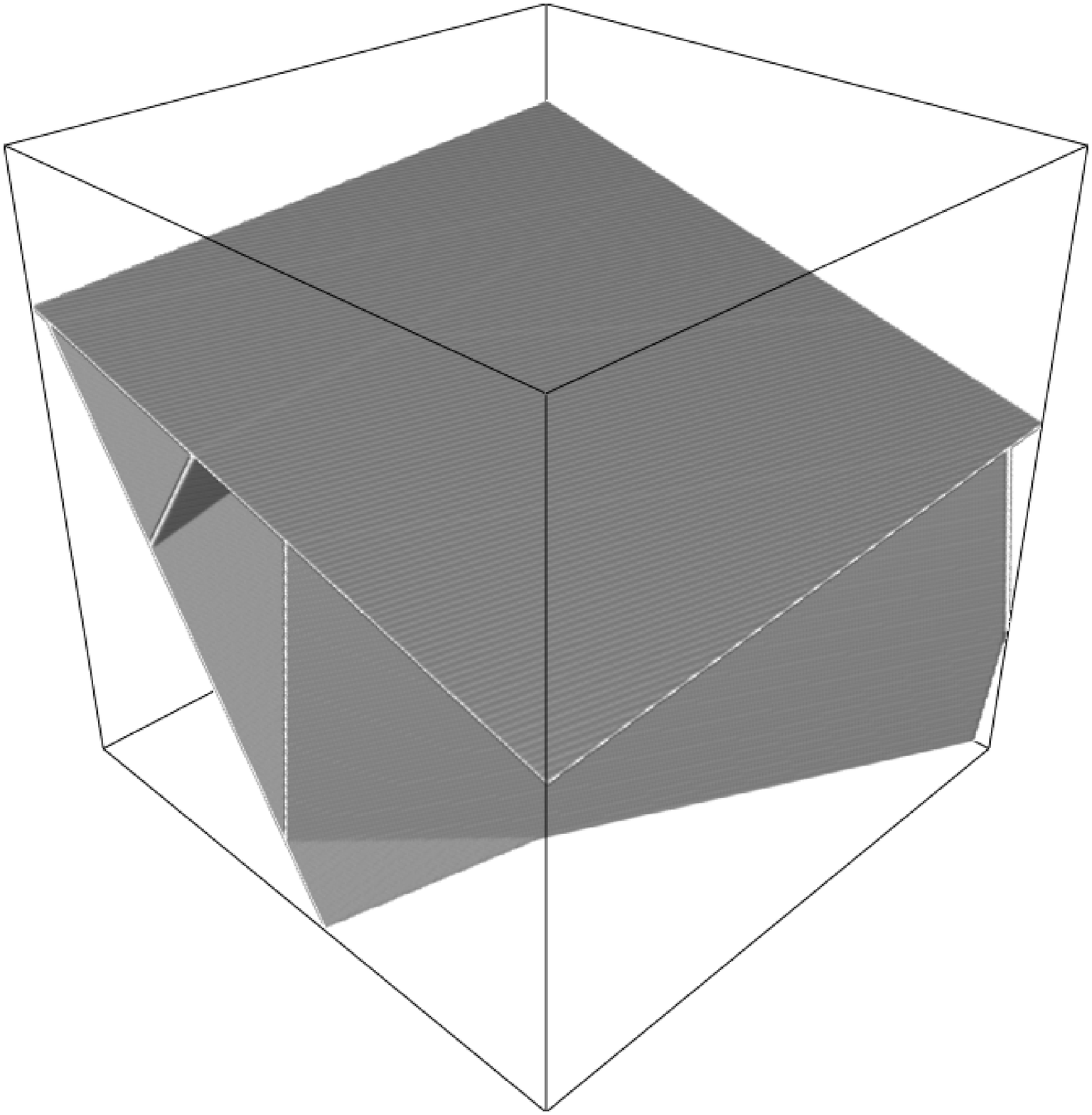}
 \includegraphics[width=0.32\columnwidth]{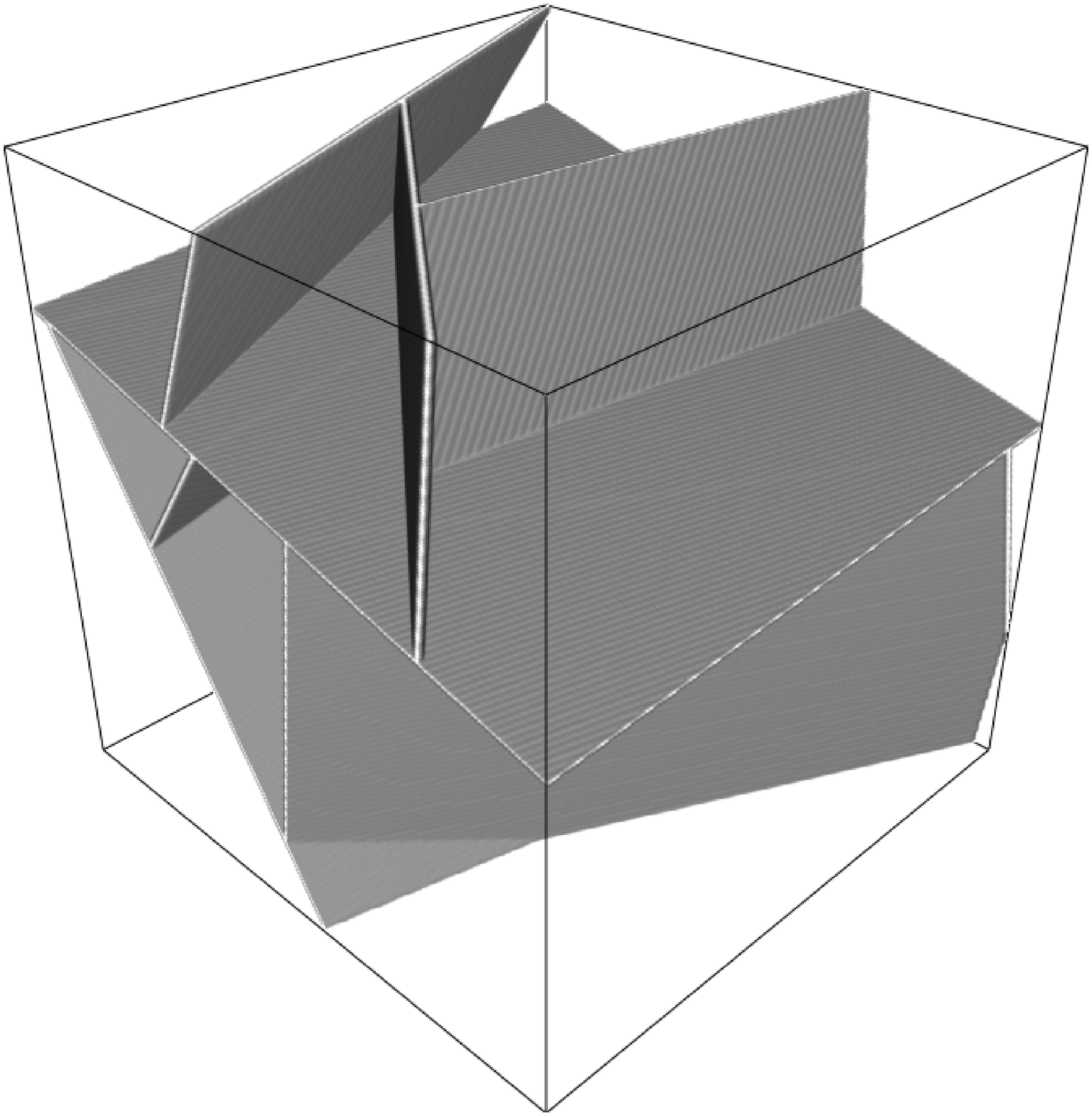}
 \caption{Illustration of the spatio-temporal construction of a STIT tessellation in a cube.}
 \label{fig:STITconstruction}
\end{center}
\end{figure}
Iteration stable tessellations (called STIT tessellations for short) have formally been introduced by Nagel and Wei\ss\ in \cite{NW05}. We confine ourself here to a less formal local description of STIT tessellations and refer to \cite{NW05} for more details. To this end, let $W\subset\RR^3$ be a compact convex polytope, and let $\L$ be a translation invariant measure on the space $\cH$ of planes in $\RR^3$ with the property that $\L([W])=1$, where $[W]=\{H\in\cH:H\cap W\neq\emptyset\}$ stands for the set of planes that hit the window $W$.

\begin{figure}
\begin{center}
 \includegraphics[width=0.5\columnwidth]{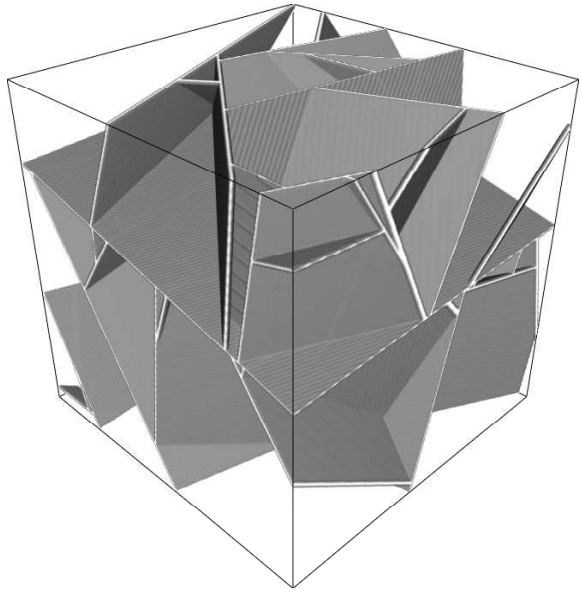}
 \caption{A spatial STIT tessellation in a cube.}
 \label{fig:STITconstructionfinal}
\end{center}
\end{figure}

We assign now to $W$ a random lifetime which is exponentially distributed with parameter $\L([W])=1$. When this time has run out, we choose a plane $H$ according to the probability distribution $\L(\cdot\cap[W])$ and divide $W$ into the two sub-cells $W^\pm=W\cap H^\pm$, where $H^\pm$ are the two half-spaces generated by $H$. These two random polytopes $W^\pm$ are now independently of each other equipped with exponentially distributed random lifetimes that have parameter $\L([W^\pm])$ and upon expiry of their lifetimes, we choose planes according to $\L(\cdot\cap [W^\pm])/\L([W^\pm])$, respectively, and subdivide $W^+$ and $W^-$ further. This construction continues independently and recursively until a previously fixed deterministic time threshold $t\in(0,\infty)$ is reached, see Figures \ref{fig:STITconstruction} and \ref{fig:STITconstructionfinal} for illustrations.

The outcome of this spatio-temporal construction, which is denoted by $Y(t,W)$, is an almost surely finite collection of random convex polytopes that subdivides the convex polytope $W$ and can be considered as a tessellation of $W$. It is a main result of \cite{NW05} that the described construction is spatially consistent. This is to say, there exists a whole space random tessellation $Y(t)$ with the property that $Y(t)\cap W$ has the same distribution as the previously constructed tessellation $Y(t,W)$. The random tessellation $Y(t)$ is usually called a \textit{STIT tessellation}. The combinatorial structure of $Y(t)$ is in the focus of our interest in the present paper.

We remark that the mean values and probabilities calculated below neither depend on the time threshold $t$ nor on the particular choice of the plane measure $\L$. For this reason these parameters are suppressed and we speak from now on about a spatial STIT tessellation and we think of it as $Y(t)$ for some $t>0$ and some fixed measure $\L$.

\subsection{Object intensities and typical objects}\label{subse:typicalobjects}

It is a crucial observation that the STIT tessellations constructed in the previous subsection are spatially homogeneous, which is to say that the distribution of the shifted tessellation $Y(t)+x$ (for some fixed measure $\L$ as above) is the same as that of $Y(t)$ for any $x\in\RR^3$. This allows us to make sense of several \textit{object intensities}. To formalize the idea, let $\sX$ be some class of random geometric objects determined by the spatial homogeneous STIT tessellations. For example $\sX$ can be the class of cells or the class of cell-separating polygons that are born in the course of the spatio-temporal construction. We assume from now on that all objects from $\sX$ are convex polytopes of a fixed dimension $i$, $0\leq i\leq 3$ (this will always be the case below), which allows us to regard $\sX$ as a subset of the space $\cP_i$ of $i$-dimensional convex polytopes. For any $p\in\cP_i$ let $c(p)$ be the midpoint of the smallest circumscribed ball, the circumcenter. Then the limit 
\begin{equation}\label{eq:objectintensities}
\la_X=\lim_{r\rightarrow\infty}{1\over r^3}\sum_{x\in\sX}{\bf 1}\{c(x) \in C_r\}
\end{equation}
is well defined, where $C_r=rC$ stands for a cube of volume $r^3$, and ${\bf 1}\{\cdot\}$  for the indicator function which is $1$ if the statement in brackets is fulfilled and $0$ otherwise, cf. \cite{SW,SKM}. It will be assumed from now on that $0<\la_X<\infty$, which will be the case for all classes considered below. In view of (\ref{eq:objectintensities}) it is clear that $\la_X$ has an interpretation as the mean number of objects of class $\sX$ per unit volume. Moreover, the spatial homogeneity of the tessellation implies that $\la_X$ does not depend on the choice and the spatial location of $C_r$.

We are now going to explain what can be considered as the typical object of a class $\sX$ as above. More precisely, we will define the \textit{distribution} of the typical object and refer to \cite{SW,SKM} for more details. Besides $\cP_i$, $0\leq i\leq 3$, define $\cP_i^o=\{p\in\cP_i:c(p)=o\}$, where $o$ denotes the origin. The space $\cP_i^o$ can be equipped with the Borel $\sigma$-field $\cB(\cP_i^o)$ generated by the usual Hausdorff distance. We define now a distribution $\bQ_X$ on $\cP_i^o$ by
\begin{equation}\label{eq:typicalobject}
\bQ_X(B)={1\over\la_X}\lim_{r\rightarrow\infty}{1\over r^3}{\sum_{x\in\sX}{\bf 1}\{c(x) \in C_r\}{\bf 1}\{x-c(x)\in B\}},
\end{equation}
where $B\in\cB(\cP_i^o)$ and $\sX$ is an object class of $i$-dimensional polytopes. A random element in $\cP_i^o$ with distribution $\bQ_X$ is called the \textit{typical object} of class $\sX$ and $\bQ_X$ is said to be the \textit{distribution of the typical object}. Intuitively it can be considered as a uniformly selected object from $\sX$ independent of its size and shape.

Besides the object intensity $\la_X$ also certain adjacency relationships are of interest for us. To define them, let $\sY$ be another class of objects determined by a spatial STIT tessellation and call $x\in\sX$ and $y\in\sY$ \textit{adjacent} if either $x\subseteq y$ or $y\subseteq x$. For fixed $x\in\sX$, we denote by $m_Y(x)$ the number of objects from $\sY$ adjacent to $x$ and further define $$\mu_{X,Y}:=\EE_{X}[m_Y(x)]:=\begingroup\textstyle\int\endgroup m_Y(x)\;\bQ_X(dx),$$ where in the expectation $\EE_{\sX}$ the subscript $\sX$ refers to integration with respect to $\bQ_X$. We interpret $\mu_{X,Y}$ as the mean number of objects from $\sY$ that are adjacent to the typical object from $\sX$. It is important to note that the first index in $\mu_{X,Y}$ is the typical object to which the mean value refers to and that the second index indicates the class of objects to which the adjacency relationship is measured. Thus, in general $\mu_{X,Y}$ and $\mu_{Y,X}$ are entirely different numbers. For example let $\sX$ be the class of cells and $\sY$ the class of vertices of the tessellation, then $\mu_{X,Y}$ denotes the mean number of vertices of the typical cell, wheres $\mu_{Y,X}$ is the mean number of cells adjacent to the typical vertex.

In the present paper we adopt the convention that an upper case $X$ refers to the typical object of class $\sX$ -- excepting the case of indices for simplicity -- and that a lower case $x$ is some representative of $\sX$. Moreover, we will in some cases consider subclasses of $\sX$, which will be denoted by $\sX[\cdot]$, where the contents of the brackets $[\cdot]$ will be a suitable suggestive symbol that is introduced in an ad hoc manner and whose meaning should be clear from the context. For example if $\sX$ is the class of vertices, $\sX[T]$ could be the collection of vertices of a special type T.

To introduce some further notation we recall that we have assumed that the members of $\sX$ are convex polytopes having a certain dimension $i\in\{0,1,2,3\}$. Then $\sX_j$ with $0\leq j<i$ will be the class of objects which are the $j$-dimensional faces of the polytopes from $\sX$ and, as above, $X_j$ stands for the typical object of $\sX_j$. For example if $\sX$ is the class of tessellation cells (here $i=3$), then $\sX_2$ or $\sX_1$ are, respectively, the class of cell facets ($2$-dimensional faces) and ridges ($1$-dimensional faces). Furthermore, by $\partial X$ and $\overset{\;\circ}{X}$ we denote, respectively, the topological boundary and the relative interior of the typical object from class $\sX$.

\section{Basic combinatorial structure}\label{sec:basicstructure}

In this section we introduce the basic objects that are in the focus of our interest and recall some basic mean values that are related to intensities and adjacencies. These are taken from \cite{NW08,TW,WC}.

\subsection{Primitive elements}\label{subsec:primitiveelements}

The primitive objects of a spatial STIT tessellation are its building blocks in some sense. By $\sZ$ (from the German word `Zellen') we denote the class of cells of the tessellation, which are the primitive objects of dimension $3$. For a point $x\in\RR^3$ let $D(x)$ be the intersection of all cells containing $x$, i.e. $\displaystyle D(x)=\bigcap_{z\in\sZ,x\in z}z$. We observe that $D(x)$ is a finite intersection of random convex polytopes, whence $D(x)$ itself is a random convex polytope. We can now introduce three other classes of primitive objects, namely $$\sV:=\{D(x):{\rm dim}\;D(x)=0,\; x\in\RR^3\},\;\sE:=\{D(x):{\rm dim}\;D(x)=1,\; x\in\RR^3\},$$ the class of \textit{vertices} and \textit{edges}, and $$\sP:=\{D(x):{\rm dim}\;D(x)=2,\; x\in\RR^3\},$$ the class of \textit{plates}. It is important to note that the primitive elements as introduced above cannot have an additional interior structure. For example, an edge or a plate cannot have any vertex in its relative interior. Nevertheless, the boundary structure of plates or cells can be rather complicated. For example a plate can be a triangle but with additional vertices on its boundary that are no corners of the triangle. The same also holds for the cells, whereby a cell can also have additional vertices and edges nested into its facets.

The relative intensities of these classes of primitive elements for a spatial STIT tessellation are $${\la_E\over\la_V}=2,\qquad{\la_P\over\la_V}={7\over 6},\qquad{\rm and}\qquad{\la_Z\over\la_V}={1\over 6}.$$ Moreover we have the following basic adjacency relationships for the primitive elements.
\begin{center}
\begin{tabular}{|cc||c|c|c|c|}
\hline
\parbox[0pt][2em][c]{0cm}{} & $\mu_{X,Y}$ & $X=V$ & $X=E$ & $X=P$ & $X=Z$\\
\hline
\parbox[0pt][2em][c]{0cm}{} & $Y=V$ & $1$ & $2$ & $36\over 7$ & $24$\\
\parbox[0pt][2em][c]{0cm}{} & $Y=E$ & $4$ & $1$ & $36\over 7$ & $36$\\
\parbox[0pt][2em][c]{0cm}{} & $Y=P$ & $6$ & $3$ & $1$ & $14$\\
\parbox[0pt][2em][c]{0cm}{} & $Y=Z$ & $4$ & $3$ & $2$ & $1$\\
\hline
\end{tabular}
\end{center}

\subsection{Facets, ridges and $P_1$-segments}\label{subsec:facetsridges}

In this subsection we consider the lower-dimensional faces of primitive elements introduced above. We start with the faces of cells. The class of cell \textit{facets} ($2$-dimensional faces) is denoted by $\sZ_2$ according to our convention and the class of \textit{ridges} ($1$-dimensional cell faces) is $\sZ_1$. Sometimes we will also refer to the elements from $\sZ_1$ as \textit{$Z_1$-segments} for simplicity. Note, moreover, that up to multiple counting the class $\sZ_0$ is the same as $\sV$. We can also consider the class $\sP_1$ of $1$-dimensional faces of tessellation plates, which are also referred to as \textit{$P_1$-segments} or sides. Finally, we notice that $\sP_0$ and $\sE_0$ are the same as $\sV$ up to multiplicities. The relative object intensities of $\sZ_1$, $\sZ_2$ and $\sP_1$ are given by $${\la_{Z_1}\over\la_V}=2,\qquad{\la_{Z_2}\over\la_V}=1\qquad{\rm and}\qquad{\la_{P_1}\over\la_V}={14\over 3}.$$ The basic adjacencies in this context are summarized in the table below (those marked with a ``$-$'' are unknown).
\begin{center}
\begin{tabular}{|cc||c|c|c|}
\hline
\parbox[0pt][2em][c]{0cm}{} & $\mu_{X,Y}$ & $X=P_1$ & $X=Z_1$ & $X=Z_2$\\
\hline
\parbox[0pt][2em][c]{0cm}{} & $Y=V$ & ${16\over 7}$ & $3$ & ${26\over 3}$\\
\parbox[0pt][2em][c]{0cm}{} & $Y=E$ & ${9\over 7}$ & $2$ & $10$\\
\parbox[0pt][2em][c]{0cm}{} & $Y=P$ & $-$ & $-$ & $7\over 3$\\
\hline
\end{tabular}
\end{center}

\subsection{I-polygons and I-segments}\label{subsec:ipolygonsandsegments}

The last interesting class of objects we introduce in this paper is the class $\sI$ of I-polygons. In the language of the spatio-temporal construction of STIT tessellations, the \textit{I-polygons} are the cell separating facets that are introduced/born during the cell dividing process. By $\sI_1$ we denote the class of \textit{I-segments}, which are the $1$-dimensional faces of I-polygons. I-polygons are also characterized by saying that they are the maximal unions of connected and coplanar plates. In the same spirit, I-segments are the maximal unions of connected and collinear tessellation edges. For the relative intensities we have $${\la_{I_1}\over\la_V}={2\over 3}\qquad{\rm and}\qquad{\la_I\over\la_V}={1\over 6}$$ and some basic adjacencies are summarized below (again, the one marked with a ``$-$'' is unknown).
\begin{center}
\begin{tabular}{|cc||c|c|}
\hline
\parbox[0pt][2em][c]{0cm}{} & $\mu_{X,Y}$ & $X=I_1$ & $X=I$\\
\hline
\parbox[0pt][2em][c]{0cm}{} & $Y=V$ & $4$ & $18$\\
\parbox[0pt][2em][c]{0cm}{} & $Y=E$ & $3$ & $24$\\
\parbox[0pt][2em][c]{0cm}{} & $Y=P$ & $-$ & $7$\\
\hline
\end{tabular}
\end{center}

\section{Fine combinatorial structure}\label{sec:finestructure}

We refine in this section the results from above by exploring several edge classifications that are based on delicate distributional results for the typical I-segment. The probabilities that the typical tessellation edge belongs to one of the introduced classes are helpful in calculating a number of new mean values, for example for the neighborhood of the typical vertex for different vertex classes.

\subsection{Vertex geometry}\label{subsec:vertex geometry}
\begin{figure}
\begin{center}
 \includegraphics[width=0.45\columnwidth]{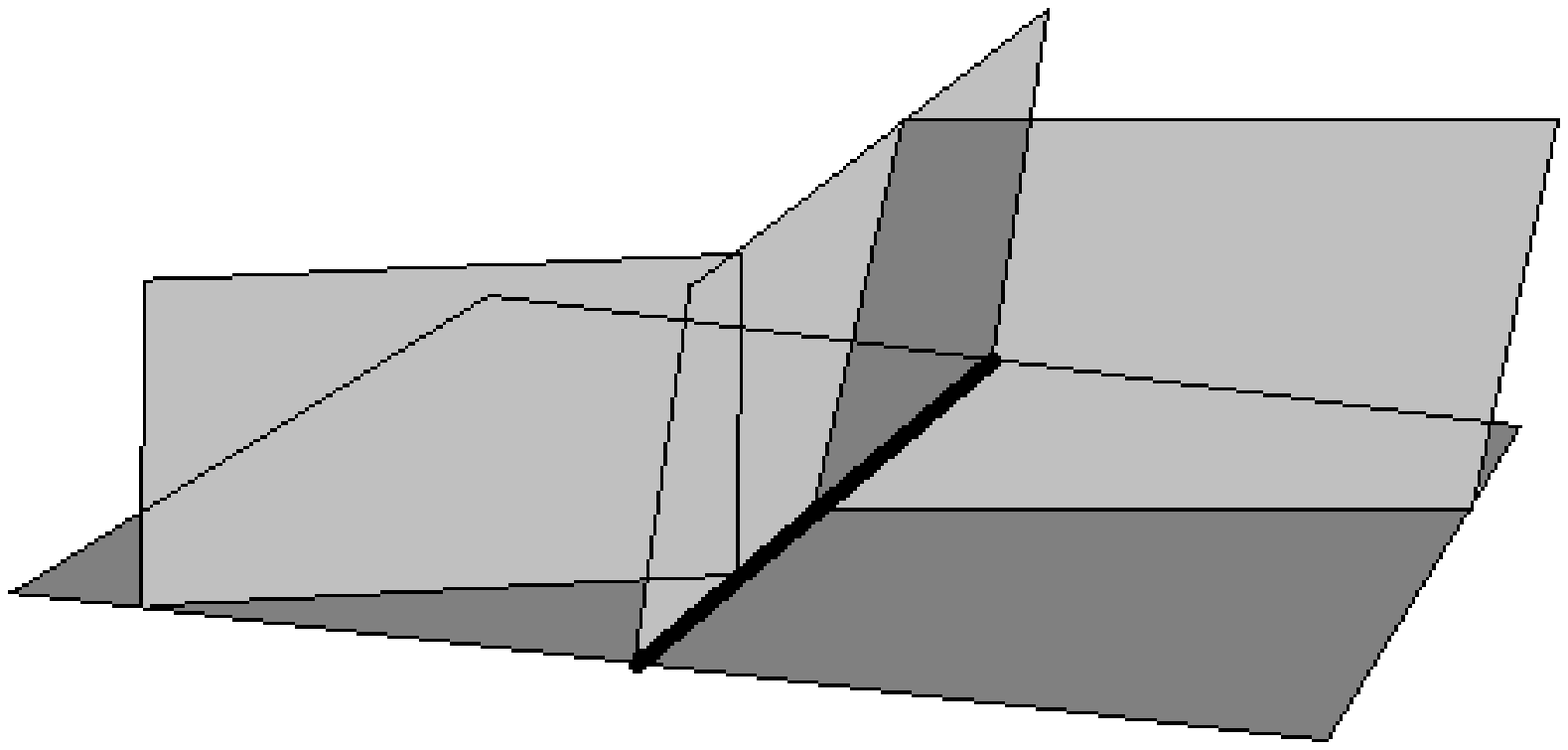}
 \includegraphics[width=0.45\columnwidth]{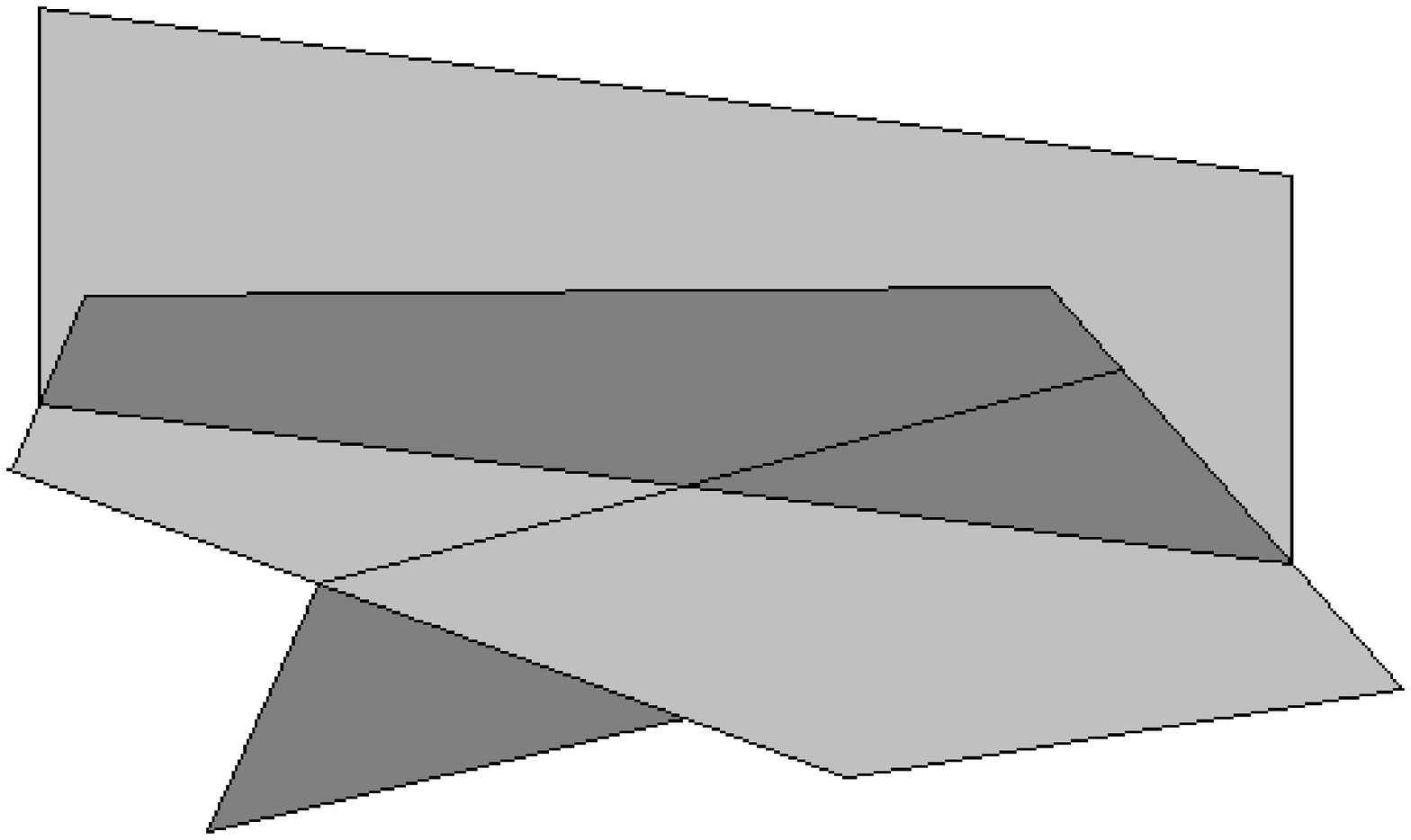}
 \caption{Two T-vertices (left) pointing to the left and to the right relative to the marked segment and a vertex of type X (right).} 
 \label{fig:vertices}
\end{center}
\end{figure}
We have introduced in Subsection \ref{subsec:primitiveelements} the class $\sV$ of vertices. In sharp contrast to the planar case, which has only vertices of type T (these are vertices in which three edges meet and exactly two of them are collinear), a spatial STIT tessellation has two types of vertices, see Figure \ref{fig:vertices}. In the course of the spatio-temporal construction from Subsection \ref{subsec:STITs} they arise as follows:
\begin{itemize}
 \item Given an already existing I-polygon, a T-vertex is generated if two further I-polygons intersect in the \textit{same} half-space determined by the I-polygon. 
 \item A vertex of type X is generated by an intersection of two further I-polygons in the two \textit{different} half-spaces specified by the given I-polygon that has been born earlier in time. 
\end{itemize}
This way, we split $\sV$ into the two subclasses $\sV=\sV[T]\cup\sV[X]$ of T- and X-vertices and denote by $$\eps_{V[T]}:=\bP(V\in\sV[T])\quad{\rm and}\quad\eps_{V[X]}:=\bP(V\in\sV[X])$$ the probability that the typical vertex is of type T or X, respectively. Denoting by $\la_{V[T]}$ and $\la_{V[X]}$ the intensity of T- and X-vertices, clearly we have  $\eps_{V[T]}={\la_{V[T]}\over\la_V}$ and $\eps_{V[X]}={\la_{V[X]}\over\la_V}$, which in view of the results from \cite{TW,WC} 
leads to 
\begin{equation}\label{eqepsvtvx}
\eps_{V[T]}={2\over 3}\qquad{\rm and}\qquad\eps_{V[X]}={1\over 3}.
\end{equation}
This may also be rephrased by saying that, on average, in a spatial STIT tessellation we have twice as many vertices of T-type compared with the X-vertices.

\subsection{Some distributional results for the typical I-segment with internal vertices}\label{subsec:isegmentdistributions1}

Our refined analysis of the combinatorial structure of spatial STIT tessellations rests on distributional properties of the typical I-segment, which have been in investigated in \cite{TWN} and will be further developed in the present subsection.

We consider the subclass $\sI_1[{\rm int}V]$ of  $\sI_1$ containing all I-segments with internal vertices and we call the typical object of $\sI_1[{\rm int}V]$ {\it the typical I-segment with internal vertices}. The breakdown of the class of vertices from Subsection \ref{subsec:vertex geometry} implies that also the vertices in the relative interior of the typical I-segment (if there are any) can either be of type T or X. Moreover, we observe that a T-vertex can be created from the left (type L) and from the right (type R) relative to the segment direction, see Figure \ref{fig:vertices} (left). So we can have three different types of vertices in the relative interior of an I-segment. We first calculate some probabilities for a randomly chosen vertex in the relative interior of the typical I-segment with internal vertices.
\begin{theorem}\label{thm:vertexinisegment}
Consider a spatial STIT tessellation.
\begin{description}
 \item[(a)] A randomly chosen vertex in the relative interior of the typical I-segment with internal vertices is of T-type with probability $p_T$ and  $$p_T={4(21+56\ln 2-54\ln 3)\over 68+208\ln 2-189\ln 3}\approx 0.433053$$ and of type X with probability $p_X$ and  $$p_X={27\ln 3-16\ln 2-16\over 68+208\ln 2-189\ln 3}\approx 0.566947.$$
 \item[(b)] A randomly chosen vertex in the relative interior of the typical I-segment with internal vertices is of type L or type R, respectively, with probability $$ p_L=p_R = \frac{1}{2}\,p_T\approx 0.216527.$$ 
\end{description}
\end{theorem}

\begin{remark}
At first sight the probabilities from part (a) of the previous theorem seem to be rather surprising, because for the mean number of internal vertices of type T or X in the typical I-segment we know $\mu_{\overset{\circ}{I_1},V[T]}=\mu_{\overset{\circ}{I_1},V[X]}=1$ from \cite{TWN}. However, the condition `having internal vertices' has an important impact on the values of $p_T$ and $p_X$.
\end{remark}

\begin{remark}
Part (b) of the previous theorem follows from the construction of STIT tessellations. The detailed argumentation given for the planar case in \cite{Cowan11} can be adopted also to our spatial setting. If a randomly chosen vertex in the relative interior of the typical I-segment with internal vertices is of T-type, then this vertex is of type L or of type R with the same probability, $p_{L|T}=p_{R|T}={1\over 2}$. In Section \ref{subsec:proofth1} we confirm this result with another short and more formal proof.
\end{remark}

Theorem \ref{thm:vertexinisegment} deals with the typical I-segment that has at least one internal vertex. Interestingly, if we condition the typical I-segment on having exactly $n\geq 1$ internal vertices, the probability for a randomly chosen vertex is of T- or X-type is no more given by $p_T$ or $p_X$ as above, but also depends on $n$. We denote it by $p_{T|n}$ and $p_{X|n}$, respectively, $n\geq 1$. This is a special feature of spatial STIT tessellations and is in sharp contrast to the much simpler planar situation as considered in \cite{Cowan11}. To calculate $p_{T|n}$ and $p_{X|n}$ we recall at first that the probability $p_n=\bP(\nu=n)$ ($n\in\NN$) that the typical I-segment has exactly $n$ vertices in its relative interior equals
\begin{equation}\label{eq:pm}
p_n=3\int\limits_0^1\int\limits_0^1(1-a)^3{(3-(1-a)(3-b))^n\over(3-(1-a)(2-b))^{n+1}}dbda,
\end{equation}
a formula that has been established in \cite{TWN}. In the present paper we will always denote by $\nu$ the number of vertices in the relative interior of the typical I-segment and also introduce the decomposition $\nu=\nu_T+\nu_X$, which takes into account the two different types of vertices. Using $p_n$ and the abbreviation $p_{m, j}=\bP(\nu_T=m,\nu_X=j)$ (see (\ref{eq:pXT}) below for an exact expression) we can formally introduce
\begin{equation}\label{eq:ptm}
p_{T|n}={1\over p_n}\sum_{k=0}^n{k\over n}p_{k,n-k}\quad{\rm and}\quad p_{X|n}=1-p_{T|n}.
\end{equation}
For example $p_{T|1}\approx 0.3854$, $p_{T|2}\approx 0.4114$ or $p_{T|20}\approx 0.6058$. This indicates that in I-segments with a few number of internal vertices it is much more likely to choose a vertex of type X compared with I-segments having a large number of vertices in their relative interior. We also have the limit relations $$\lim_{n\rightarrow\infty} p_{T|n}={2\over 3}\quad{\rm and}\quad\lim_{n\rightarrow\infty} p_{X|n}={1\over 3},$$ which are exactly the probabilities $\eps_{V[T]}$ and $\eps_{V[X]}$ from (\ref{eqepsvtvx}). An analytic expression for $p_{T|n}$ is provided by $$p_{T|n}={6\over np_n}\left[{F_1(2;n,-n;3;-1/2,-2)\over 2^{n+1}}-{3F_1(3;n,-n;4;-1/2,-2)\over 2^n}\right.$$ $$\qquad\qquad\left.+\;3^n\left({_2F_1(n,n+3;n+4;-2)\over n+3}-{_2F_1(n,n+2;n+3;-2)\over n+2}\right)\right],$$ where $F_1$ and $_2F_1$ are, respectively, Appell's and Gau\ss's hypergeometric functions. We will see in the next subsection that the numbers $p_{T|n}$ and $p_n$ will play an important r\^ole in our fine combinatorial analysis.

It is of great importance that the probability for choosing a vertex of T- or X-type on the typical I-segment with an fixed number of internal vertices depends on that number of internal vertices. Another crucial observation is that -- by construction -- the types of two neighbouring internal vertices are independent and identically distributed. For example, the probability that on the typical I-segment with $\nu=n\geq 2$ internal vertices we find two neighboring T-vertices is $p_{T|n}^2$. Thus, looking along  a typical I-segment, the vertex type is determined by independent coin tosses with probabilities $p_{T|n}$ and $p_{X|n}=1-p_{T|n}$. Moreover, if a vertex is of type T, then the labels L (left) and R (right) are determined by independent coin tosses with a fair coin. This will be used frequently later in this paper. (Formally this can be established by standard properties of Poisson point processes, which arise as intersection of a spatial STIT tessellation with a line, see \cite{Cowan11} for a similar argument in the planar case.)

\subsection{Classification of edges}\label{subsec:edgeclassification}

We now refine the analysis initiated in Subsection \ref{subsec:primitiveelements} by dividing the class $\sE$ of edges into subclasses, where we explore several possibilities. The first one is based on a classification according to the type of endpoints, whereas the other two classifications concentrate on the equality relationship with other geometric objects. It will be shown that the probability that the typical edge belongs to one of the considered subclasses of $\sE$ can be expressed in terms of $p_n$ and $p_{T|n}$ from (\ref{eq:pm}) and (\ref{eq:ptm}). It is also worth noticing that the probability that the typical edge belongs to one of the subclasses of $\sE$ is the same as the proportion of tessellation edges that belong to this subclass.

We would like to point out that in \cite{Cowan11} a similar study has been carried out for the planar STIT tessellation, where we have only one type of vertices. Now, the situation with two different types of vertices in the spatial case is much more tricky, produces a lot of new effects and allows us to consider several different classifications.

\paragraph{Classification according to endpoints.} The first way is to classify an edge by the type of its two endpoints. Recall that a spatial STIT tessellation has T- and X-type vertices and so an edge can have two T-vertices, two X-vertices or one T- and one X-vertex as its endpoints. This way, the class $\sE$ of edges can be split into the three subclasses $\sE[TT]$, $\sE[XX]$ and $\sE[TX]$. We are interested in the probabilities that the typical edge belongs to one of these classes. More formally, we define $$\eps_{E[TT]}:=\bP(E\in\sE[TT])=\frac{\lambda_{E[TT]}}{\lambda_{E}},\quad \eps_{E[XX]}:=\bP(E\in\sE[XX])=\frac{\lambda_{E[XX]}}{\lambda_{E}}$$ and $$\eps_{E[TX]}:=\bP(E\in\sE[TX])=\frac{\lambda_{E[TX]}}{\lambda_{E}}$$ which clearly satisfy $\eps_{E[TT]}+\eps_{E[XX]}+\eps_{E[TX]}=1$.
\begin{theorem}\label{thm:TXclassification}
For any spatial STIT tessellation we have
\begin{eqnarray}
\nonumber \eps_{E[TT]} &=& {1\over 3}p_0+{2\over 3}p_{T|1}p_1+{1\over 3}\sum_{n=2}^\infty\left(2p_{T|n}+(n-1)p_{T|n}^2\right)p_n\\
\nonumber & & \approx 0.442878,\\
\nonumber \eps_{E[XX]} &=& \eps_{E[TT]}-{1\over 3}\approx 0.109545,\\
\nonumber \eps_{E[TX]} &=& {4\over 3}-2\eps_{E[TT]}\approx 0.447577,
\end{eqnarray}
where $p_n$ is given by (\ref{eq:pm}) and $p_{T|n}$  by (\ref{eq:ptm}).
\end{theorem}

Unfortunately, we were not able to compute the exact value of $\eps_{E[TT]}$, which is mainly due to the fact that the probabilities $p_{T|n}$ and $p_n$ have rather complicated expressions and, moreover, that the squared value of $p_{T|n}$ and even $1/p_n$ enters the formula. A similar remark also applies to the probabilities in the Theorems \ref{thm:P1classification} and \ref{thm:Z1classification} below.

To obtain our numerical values we have used Mathematica and have evaluated the first $10\;000$ terms of the sums by the implemented numerical integration methods. To truncate the sums for considerably smaller $n$, $n=100$ say, does not lead to accurate results, because of the very large tails of the distribution $\{p_n\}_{n\in\NN}$.
\begin{figure}
\begin{center}
 \includegraphics[width=5cm]{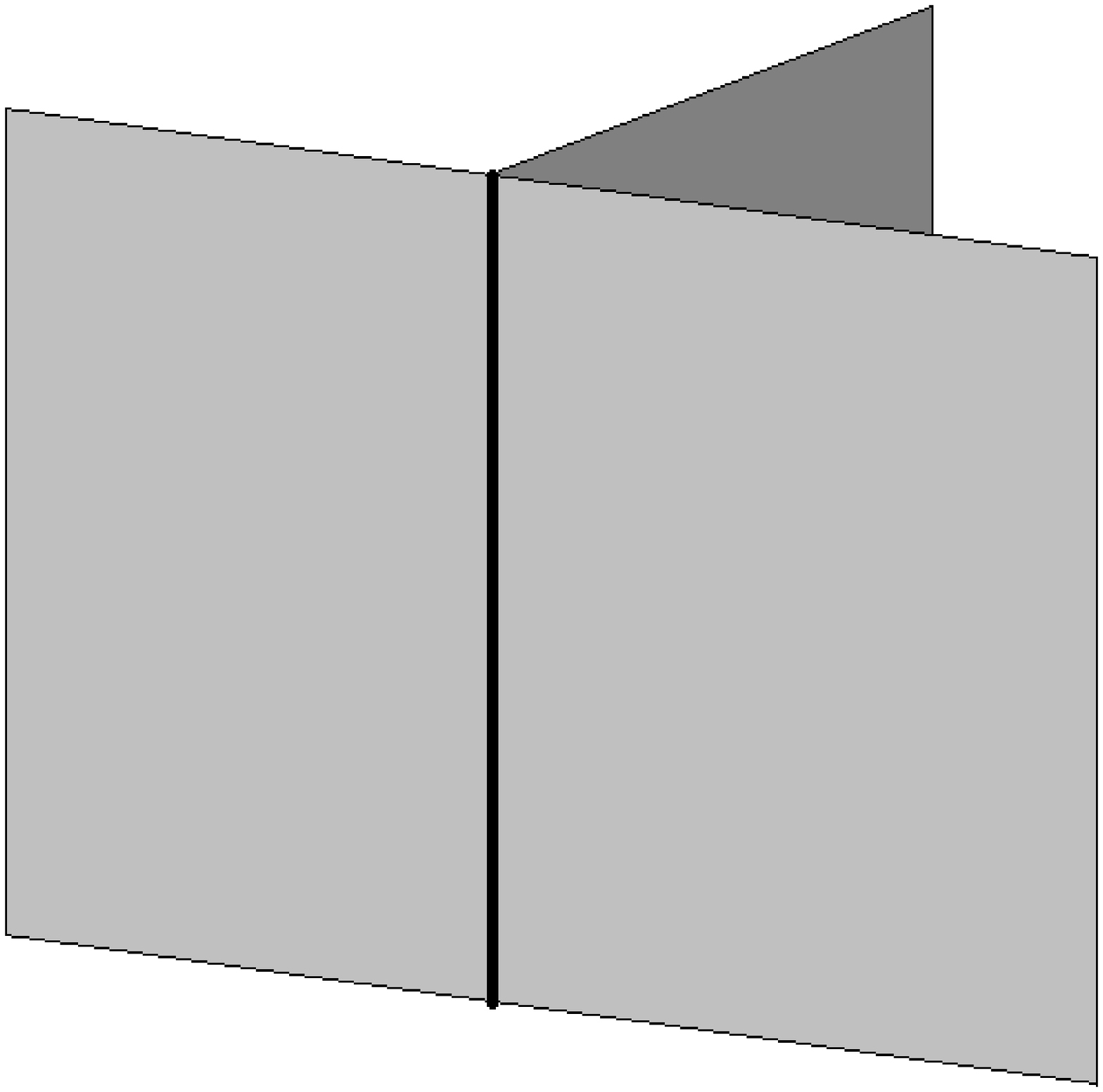}\qquad\qquad
 \includegraphics[width=5.5cm]{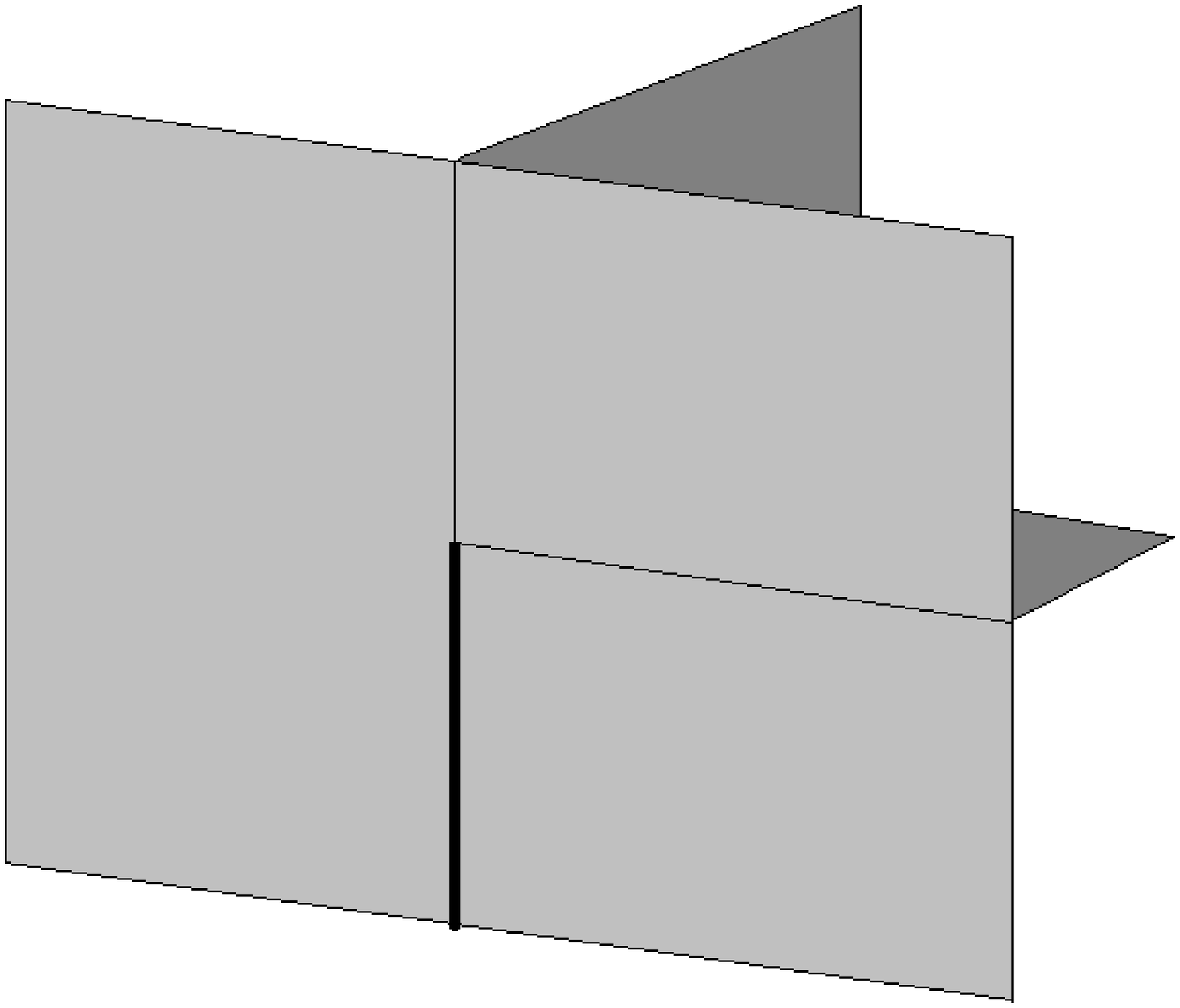}\\
 \includegraphics[width=5.8cm]{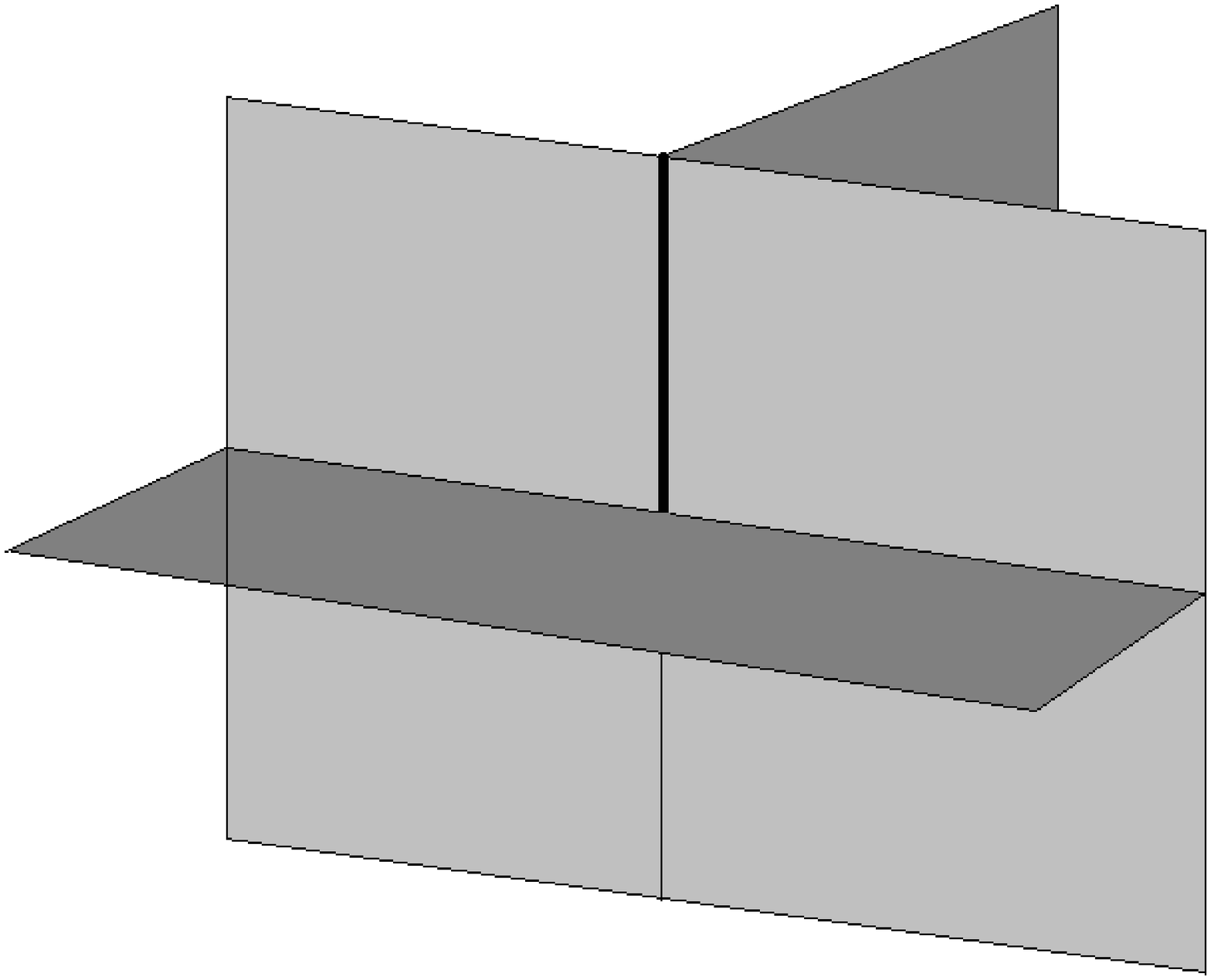}
 \includegraphics[width=5.8cm]{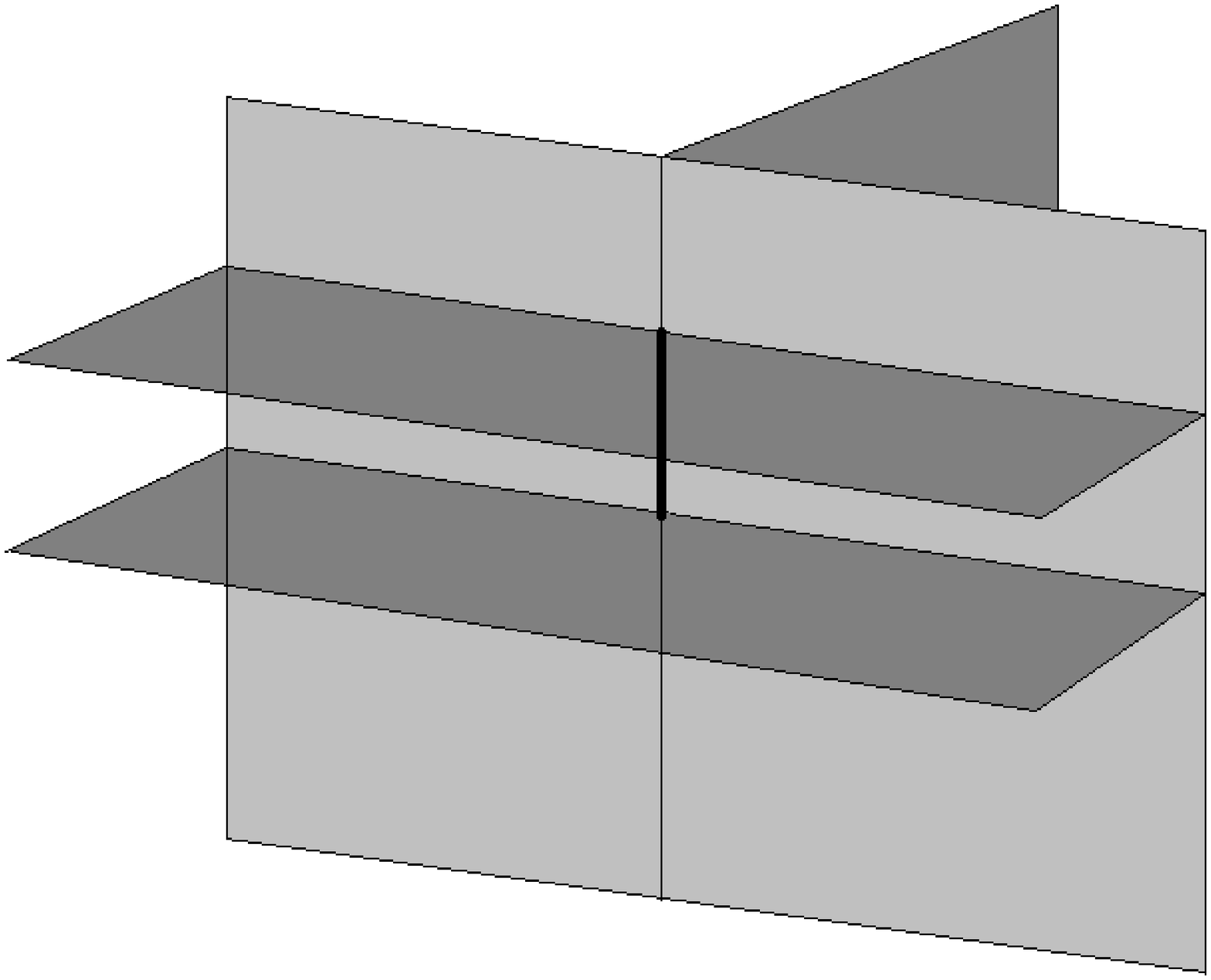}\\
 \includegraphics[width=6.2cm]{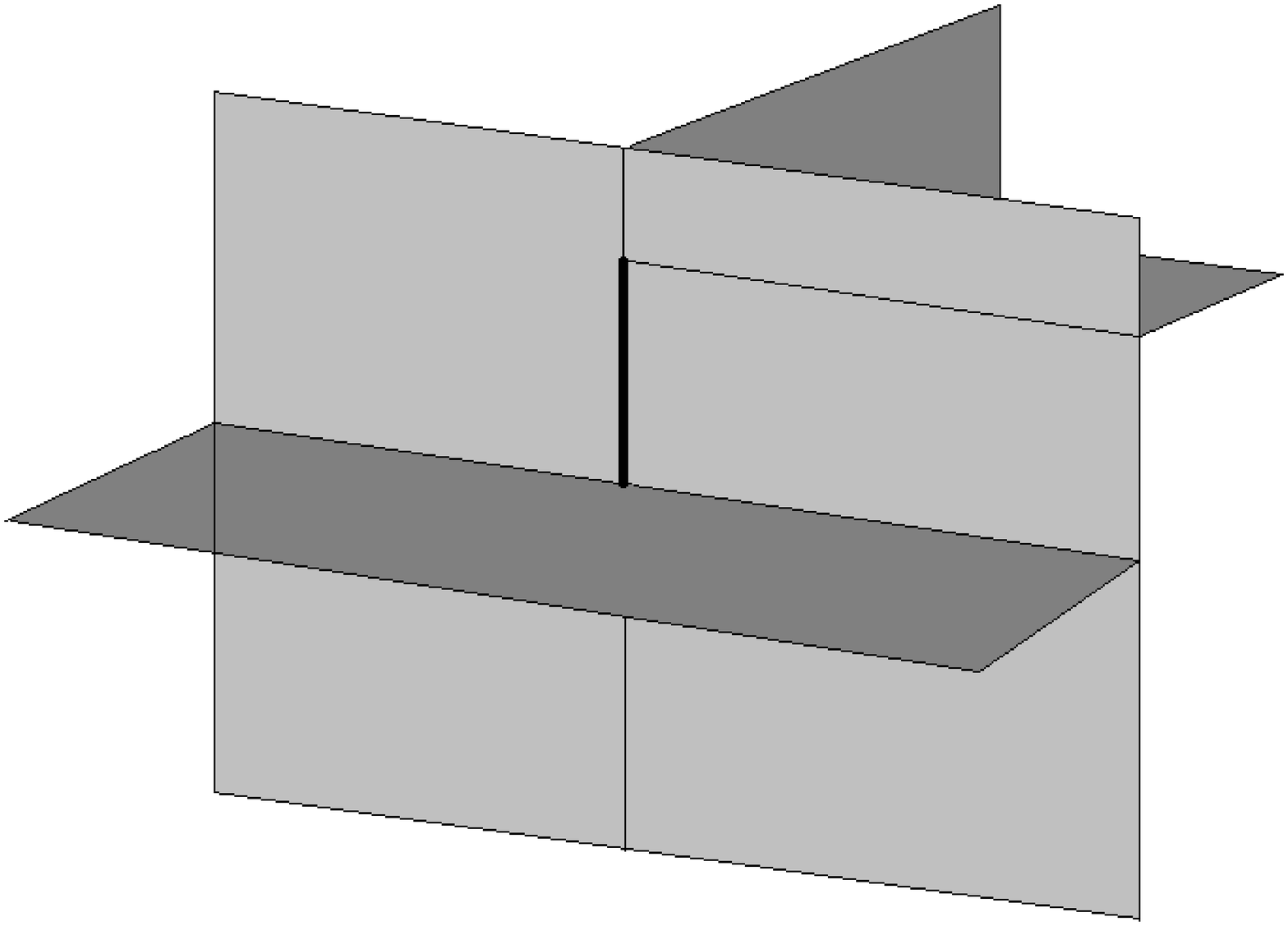}
 \includegraphics[width=5.5cm]{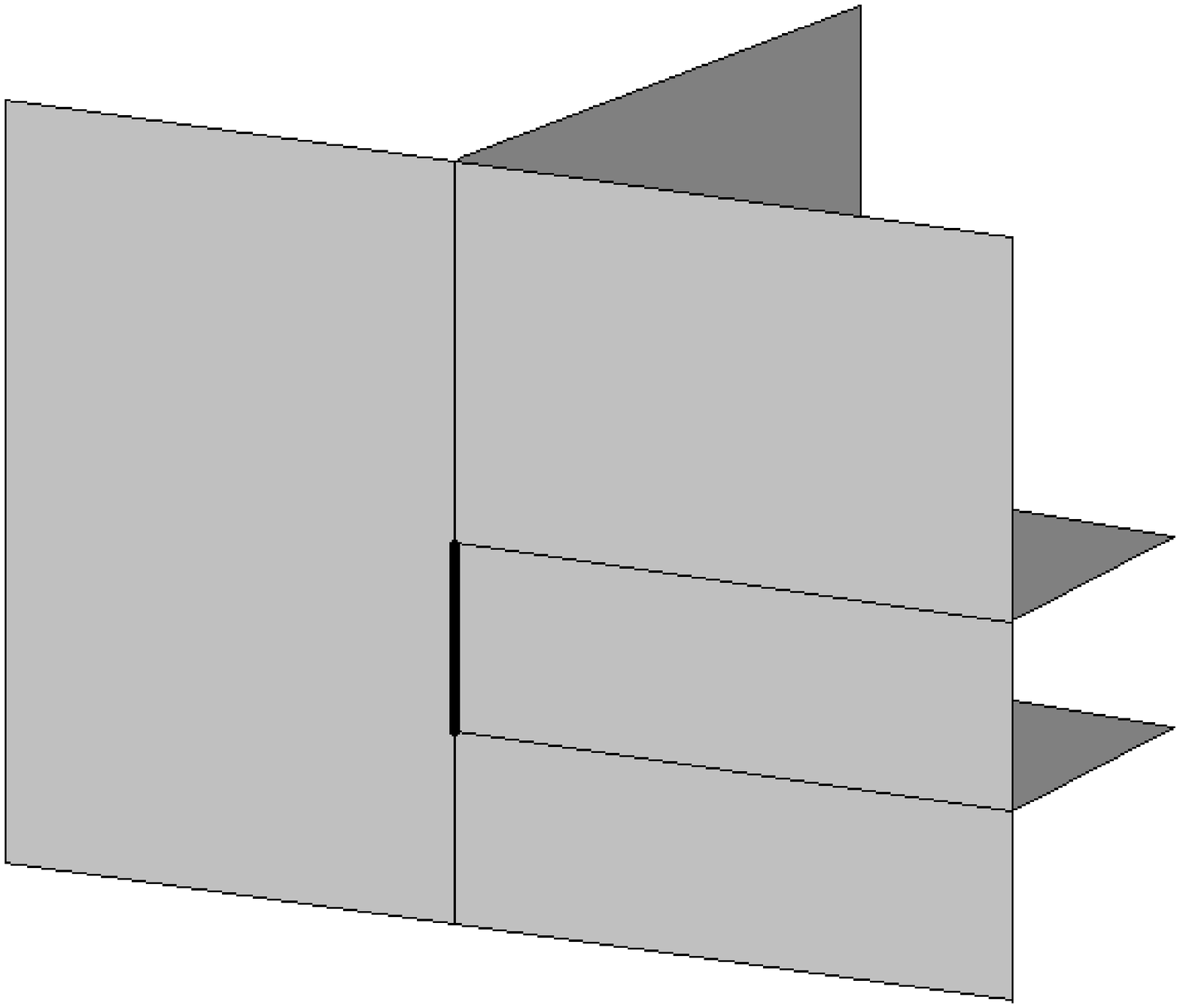}\\
 \includegraphics[width=5.5cm]{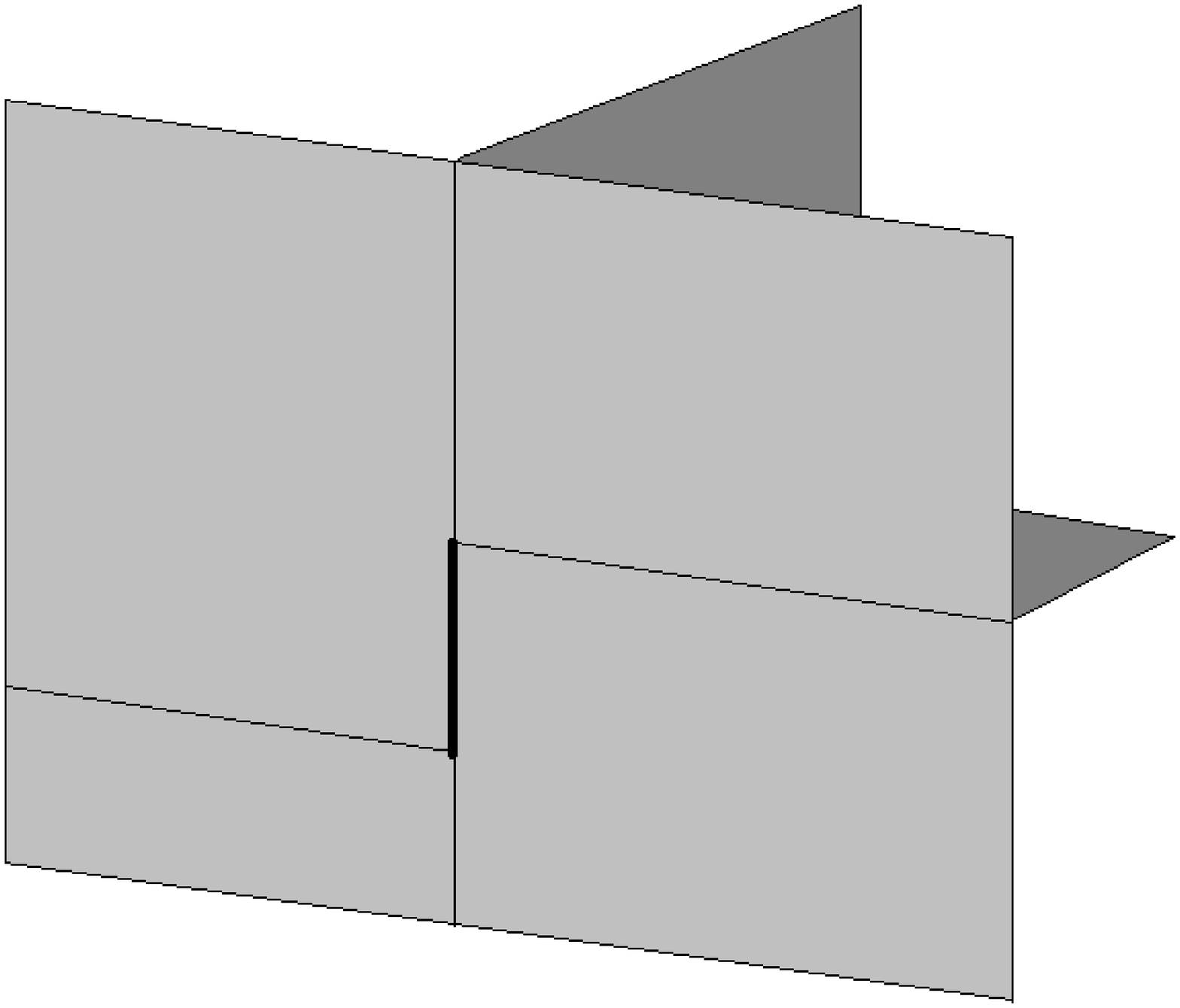}
 \caption{Illustration of different edge types on an I-segment of a spatial STIT tessellation. First row (left): $E[TT]$, $E[P_1,3]$,$E[Z_1,2]$; first row (right): $E[TT]$, $E[P_1,2]$, $E[Z_1,1]$; second row (left): $E[TX]$, $E[P_1,2]$, $E[Z_1,0]$; second row (right): $E[XX]$, $E[P_1,2]$, $E[Z_1,0]$; third row (left): $E[TX]$, $E[P_1,1]$, $E[Z_1,0]$; third row (right): $E[TT]$, $E[P_1,2]$, $E[Z_1,1]$; fourth row: $E[TT]$, $E[P_1,1]$, $E[Z_1,0]$.} 
 \label{fig:edges}
\end{center}
\end{figure}

\paragraph{Classification according to equality with $P_1$-segments.} Any edge of a STIT tessellation is adjacent to three plates. Our next edge classification is based on the observation that an edge can be \textit{equal} (not adjacent) to either one, two or three plate sides, the $P_1$-segments, whereas it is not possible, that an edge is not equal to any plate side, see Figure \ref{fig:edges}. In this figure we illustrate the different edge types on an I-segment. We start with an I-segment having no internal vertices, then with an I-segment that has one internal vertex of type T or X, respectively, etc. In the classification according to equality with $P_1$-segments we divide the class $\sE$ into the subclasses $\sE[P_1,1]$, $\sE[P_1,2]$ and $\sE[P_1,3]$ and define $$\eps_{E[P_1,i]}:=\bP(E\in\sE[P_1,i]),\quad i=1,2,3$$ as the probability that the typical edge belongs to class $\sE[P_1,i]$, where $i$ stands for the number of plate sides equal to the edge. Again we can say that $\eps_{E[P_1,i]}$ denotes the proportion of edges of the STIT tessellation which are equal to $i$ plate sides.
\begin{theorem}\label{thm:P1classification}
A spatial STIT tessellation satisfies
\begin{eqnarray}
\nonumber \eps_{E[P_1,1]} &=& {2\over 3}p_1+{1\over 3}\sum_{n=2}^\infty\left(2+(n-1)p_{T|n}^2/2+(n-1)p_{X|n}^2\right)p_n\\
\nonumber & & \approx 0.555046,\\
\nonumber \eps_{E[P_1,2]} &=& {1\over 3}\sum_{n=2}^\infty\left((n-1)p_{T|n}^2/2+2(n-1)p_{T|n}p_{X|n}\right)p_n\approx 0.300657,\\
\nonumber \eps_{E[P_1,3]} &=& {1\over 3}p_0={68\over 8}\ln 3-{26\over 3}\ln 2-{5\over 2} \approx 0.144296
\end{eqnarray}
with  $p_n$ and $p_{T|n}$ given, respectively, by (\ref{eq:pm}) and (\ref{eq:ptm}).
\end{theorem}

\paragraph{Classification according to equality with ridges.} Another breakdown of $\sE$ can be based on the equality relationship with ridges ($Z_1$-segments) instead of $P_1$-segments. In this context we observe that an edge is adjacent to two ridges of the spatial STIT tessellation, so it can be equal to either two, one or no ridge, see Figure \ref{fig:edges}. We can this way decompose $\sE$ into $\sE=\sE[Z_1,0]\cup\sE[Z_1,1]\cup\sE[Z_1,2]$ and introduce $$\eps_{E[Z_1,i]}:=\bP(E\in\sE[Z_1,i]),\quad i=0,1,2,$$ where $\sE[Z_1,i]$ is the class of tessellation edges that are equal to exactly $i$ ridges ($i=0,1,2$).
\begin{theorem}\label{thm:Z1classification}
The values for $\eps_{E[Z_1,i]}$ ($i=0,1,2$) are given by
\begin{eqnarray}
\nonumber \eps_{E[Z_1,0]} &=& {2\over 3}p_1p_{X|1}+{1\over 3}\sum_{n=2}^\infty\left(2p_{X|n}+(n-1)p_{X|n}^2+2(n-1)p_{T|n}p_{X|n}\right.\\
\nonumber & & \left.+\;(n-1)p_{T|n}^2/2\right)p_n\approx 0.624550,\\
\nonumber \eps_{E[Z_1,1]} &=& {2\over 3}p_1p_{T|1}+{1\over 3}\sum_{n=2}^\infty\left(2p_{T|n}+(n-1)p_{T|n}^2/2\right)p_n\approx 0.231154\\
\nonumber \eps_{E[Z_1,2]} &=& {1\over 3}p_0={68\over 8}\ln 3-{26\over 3}\ln 2-{5\over 2} \approx 0.144296,
\end{eqnarray}
where again $p_n$ and $p_{T|n}$ are as in (\ref{eq:pm}) and (\ref{eq:ptm}).
\end{theorem}

It is important to note that the presented edge classifications are essentially different from each other and that any of the probabilities calculated above carries important information about the structure of a spatial STIT tessellation. We will in the next subsection use these $\eps$-values to obtain new mean geometric values that arise in a refined analysis of the combinatorial structure of the model under consideration. 

\begin{remark}
It is worth noticing that with the $\eps$-values from Theorem \ref{thm:P1classification} and \ref{thm:Z1classification} and simple mean value relations further $\eps$-values can be calculated, for example the probability that the typical plate-side ($P_1$-segment) and the typical ridge ($Z_1$-segment), respectively, is equal (not adjacent) to one edge, $\eps_{P_1[E,1]}$ and $\eps_{Z_1[E,1]}$. We obtain
$$\eps_{P_1[E,1]}=\frac{\lambda_E}{\lambda_{P_1}}(\eps_{E[P_1,1]}+2\eps_{E[P_1,2]}+3\eps_{E[P_1,3]})\approx 0.681106$$
and
$$\eps_{Z_1[E,1]}=\frac{\lambda_E}{\lambda_{Z_1}}(\eps_{E[Z_1,1]}+2\eps_{E[Z_1,2]})\approx 0.519746.$$
Another interesting interpretation of these values for the typical plate side $P_1$ and the typical ridge $Z_1$ is the following:
\begin{eqnarray}
\nonumber \eps_{P_1[E,1]} &=& \bP(P_1\;{\rm has\; no\; internal\; vertices}) ,\\
\nonumber \eps_{Z_1[E,1]} &=& \bP(Z_1\;{\rm has\; no\; internal\; vertices}).
\end{eqnarray}
\end{remark}

\subsection{New geometric mean values}\label{subsec:newmeanvalues}

In this section we calculate a number of new mean values for geometric parameters of a spatial STIT tessellation and continue the studies initiated in \cite{NW08,TW}. However, we emphasize that these mean values were not available before as they rely on the new probabilities from Theorems \ref{thm:TXclassification}--\ref{thm:Z1classification} above.

\paragraph{Vertex-edge adjacencies.} In this paragraph we are dealing with certain vertex-edge adjacencies and start with the mean values $\mu_{V,E[TT]}$, $\mu_{V,E[XX]}$ and $\mu_{V,E[TX]}$ and those for the typical T- and X-vertex.

Let us first consider the typical T-vertex. With the mean value relation $\lambda_{V[T]} \mu_{V[T],E[TT]}  = 2 \lambda_{E[TT]}$ and with the relation $\lambda_{V[T]} = \frac{2}{3} \lambda_V = \frac{1}{3} \lambda_E$ for the intensities, see Section \ref{subsec:primitiveelements} and Section \ref{subsec:vertex geometry}, we obtain for the mean number of TT-edges emanating from the typical T-vertex:
$$\mu_{V[T],E[TT]} = 6 \frac{\lambda_{E[TT]}}{\lambda_E} = 6 \eps_{E[TT]}. $$
For the other edge types it is easy to see, that
$$\mu_{V[T],E[XX]} = 0 \qquad {\rm and} \qquad \mu_{V[T],E[TX]} = 4- 6 \eps_{E[TT]}, $$
because $\mu_{V[T],E} = 4$. Similar considerations for the typical X-vertex yield
$$\mu_{V[X],E[XX]} = 2 \frac{\lambda_{E[XX]}}{\lambda_E} = 12 \eps_{E[XX]}=12 \eps_{E[TT]} -4 $$
and
$$\mu_{V[X],E[TT]} = 0, \qquad \mu_{V[X],E[TX]} = 8- 12 \eps_{E[TT]}.$$
For the typical vertex we have
$$ \lambda_V \mu_{V,E[TT]} = \lambda_{V[T]}  \mu_{V[T],E[TT]} \ \Longrightarrow \ \mu_{V,E[TT]} = 4 \eps_{E[TT]} $$
and
$$ \lambda_V \mu_{V,E[XX]} = \lambda_{V[X]}  \mu_{V[X],E[XX]} \ \Longrightarrow \ \mu_{V,E[XX]} = 4 \eps_{E[TT]}- \frac{4}{3} $$
and finally
$$ \mu_{V,E[TX]} = \frac{16}{3} - 8 \eps_{E[TT]}.$$ 
With Theorem \ref{thm:TXclassification} the numerical values are summarized in the table below.
\begin{center}
\begin{tabular}{|cc||c|c|c|}
\hline
\parbox[0pt][2em][c]{0cm}{} & & $TT$ & $XX$ & $TX$\\
\hline
\parbox[0pt][2em][c]{0cm}{} & $\mu_{V[T],E[\cdot]}$ & $2.65727$ & $0$ & $1.34273$\\
\parbox[0pt][2em][c]{0cm}{} & $\mu_{V[X],E[\cdot]}$ & $0$ & $1.31454$ & $2.68546$\\
\parbox[0pt][2em][c]{0cm}{} & $\mu_{V,E[\cdot]}$ & $1.77151$ & $0.43817$ & $1.79031$\\
\hline
\end{tabular}
\end{center}
We deal now with the mean adjacencies $\mu_{V,E[P_1,i]}$ ($i=1,2,3$) and $\mu_{V,E[Z_1,j]}$ ($j=0,1,2$) for the typical vertex.
With $\lambda_V \mu_{V,E[P_1,i]}= 2 \lambda_{E[P_1,i]}$ and a similar relation for the edges from $\sE[Z_1,j]$ we have the exact values $$\mu_{V,E[P_1,i]}=4\eps_{E[P_1,i]}\;(i=1,2,3)\quad{\rm and}\quad\mu_{V,E[Z_1,j]}=4\eps_{E[Z_1,j]}\;(j=0,1,2).$$ They lead with Theorem \ref{thm:P1classification} and \ref{thm:Z1classification} to the following numerical results.

\begin{center}
\begin{tabular}{|cc||c|c|c||c|c|c|}
\hline
\parbox[0pt][2em][c]{0cm}{} & & $E[P_1,1]$ & $E[P_1,2]$ & $E[P_1,3]$ & $E[Z_1,0]$ & $E[Z_1,1]$ & $E[Z_1,2]$\\
\hline
\parbox[0pt][2em][c]{0cm}{} & $\mu_{V,\cdot}$ & $2.22019$ & $1.20263$ & $0.57718$ & $2.49820$ & $0.92461$ & $0.57718$\\
\hline
\end{tabular}
\end{center}

The adjacencies between the typical T-vertex or the typical X-vertex, respectively, and the edges from the classes $\sE[P_1,i]$ and $\sE[Z_1,j]$ are more complicated.

All edges in the subclasses $\sE[P_1,3]$, $\sE[Z_1,1]$ and $\sE[Z_1,2]$ are TT-edges. With the mean value relations $\lambda_{V[T]} \mu_{V[T],E[P_1,3]}= 2 \lambda_{E[P_1,3]}$ and $\lambda_{E[P_1,3]}= \lambda_E \eps_{E[P_1,3]}$ and the two intensity relations $\lambda_E=2{\lambda_V}$ and $\lambda_{V[T]}=\frac{2}{3}{\lambda_V}$ we obtain for the mean number of $E[P_1,3]$-edges emanating from the typical T-vertex
$$ \mu_{V[T],E[P_1,3]} = 6 \eps_{E[P_1,3]}$$
and analogously
$$ \mu_{V[T],E[Z_1,1]} = 6 \eps_{E[Z_1,1]} \quad {\rm and} \quad \mu_{V[T],E[Z_1,2]} = 6 \eps_{E[Z_1,2]},$$
and of course $ \mu_{V[X],E[P_1,3]} =\mu_{V[X],E[Z_1,1]}=\mu_{V[X],E[Z_1,2]}=0.$
From $$\mu_{V[T],E[Z_1,0]}+\mu_{V[T],E[Z_1,1]}+\mu_{V[T],E[Z_1,2]}=\mu_{V[T],E}=4$$ we infer
$$\mu_{V[T],E[Z_1,0]} = 4 - 6 (\eps_{E[Z_1,1]}+\eps_{E[Z_1,2]}) $$
and $\mu_{V[X],E[Z_1,0]} = 4$ for the typical X-vertex. The last result is obvious from the X-vertex geometry, because an X-vertex is an internal vertex of four ridges and therefore all emanating edges can not be equal to a ridge, see also Figure \ref{fig:vertices}.
To find the relations for $ \mu_{V[\cdot],E[P_1,1]}$ and $ \mu_{V[\cdot],E[P_1,2]}$ we need the same technique as the one used in the proof of Theorems \ref{thm:TXclassification}-\ref{thm:Z1classification}. Here, we will state only the results, the proofs are postponed to Subsection \ref{subsec:proofvertexedge}. We have
\begin{eqnarray}
\nonumber \mu_{V[T],E[P_1,1]} &=& 27\ln 3-28\ln 2-{19\over 2},\\
\nonumber \mu_{V[X],E[P_1,1]} &=& 2 \sum_{n=2}^{\infty} (n-1)  p_{T|n}p_{X|n}p_n,
\end{eqnarray}
and, moreover, 
\begin{eqnarray}
\nonumber \mu_{V[T],E[P_1,2]} &=& 4-  6 \eps_{E[P_1,3]} -\sum_{n=2}^{\infty} (n-1)  p_{T|n}p_n,\\
\nonumber \mu_{V[X],E[P_1,2]} &=& 4-2 \sum_{n=2}^{\infty} (n-1)  p_{T|n}p_{X|n}p_n,
\end{eqnarray}
where $p_n$ and $p_{T|n}$ are as in (\ref{eq:pm}) and (\ref{eq:ptm}). We can now calculate the numerical values for these adjacency mean values with Theorem \ref{thm:P1classification}-\ref{thm:Z1classification}. 
\begin{center}
\begin{tabular}{|cc||c|c|c||c|c|c|}
\hline
\parbox[0pt][2em][c]{0cm}{} & & $E[P_1,1]$ & $E[P_1,2]$ & $E[P_1,3]$ & $E[Z_1,0]$ & $E[Z_1,1]$ & $E[Z_1,2]$\\
\hline
\parbox[0pt][2em][c]{0cm}{} & $\mu_{V[T],\cdot}$ & $0.75441$ & $2.37981$ & $0.86577$ & $1.74730$ & $1.38692$ & $0.86577$\\
\hline
\end{tabular}
\end{center}
\begin{center}
\begin{tabular}{|cc||c|c|c||c|c|c|}
\hline
\parbox[0pt][2em][c]{0cm}{} & & $E[P_1,1]$ & $E[P_1,2]$ & $E[P_1,3]$ & $E[Z_1,0]$ & $E[Z_1,1]$ & $E[Z_1,2]$\\
\hline
\parbox[0pt][2em][c]{0cm}{} & $\mu_{V[X],\cdot}$ & $0.69968$ & $3.30032$ & $0$ & $4$ & $0$ & $0$\\
\hline
\end{tabular}
\end{center}

\paragraph{The neighborhood of typical vertices.} Any vertex $v\in\sV$ of a STIT tessellation has $4$ adjacent edges, which have one endpoint equal to $v$ and the other one is called a \textit{neighbor} of $v$. So, each vertex has four neighboring vertices. We are interested in the mean number of T-vertices and of X-vertices in the neighborhood of the typical vertex and more specifically, also of the typical T- and the typical X-vertex. We denote these mean values by $\eta_{V,V[\cdot]}$, $\eta_{V[T],V[\cdot]}$ and $\eta_{V[X],V[\cdot]}$, respectively, which are defined in the same spirit as the $\mu$-mean values in Subsection \ref{subse:typicalobjects} and $[\cdot]$ stands for $[T]$ or $[X]$, respectively. Of course we have 
$$\eta_{V[T],V[\cdot]}= \mu_{V[T],E[T\cdot]}, \qquad \eta_{V[X],V[\cdot]}= \mu_{V[X],E[X\cdot]}$$
and similarly
$$\eta_{V,V[\cdot]}= \eps_{V[T]}\mu_{V[T],E[T\cdot]} + \eps_{V[X]}\mu_{V[X],E[X\cdot]}.$$
Using the results from the vertex-edge adjacencies we obtain for the typical vertex
$$\eta_{V,V[T]}=8/3 \qquad {\rm and} \qquad \eta_{V,V[X]}=4/3$$
and for the typical T- and X-vertex
 \begin{alignat}{4}
\nonumber &\eta_{V[T],V[T]} && = 6\eps_{E[TT]},\qquad &&\eta_{V[T],V[X]} &&= 4-6\eps_{E[TT]},\\
\nonumber &\eta_{V[X],V[T]} && = 8-12\eps_{E[TT]},\qquad &&\eta_{V[X],V[X]} &&= 12\eps_{E[TT]}-4.
\end{alignat}
With Theorem \ref{thm:TXclassification} we find the following values.
\begin{center}
\begin{tabular}{|cc||c|c|c|}
\hline
\parbox[0pt][2em][c]{0cm}{} & & $T$ & $X$\\
\hline
\parbox[0pt][2em][c]{0cm}{} & $\eta_{V[T],V[\cdot]}$ & $2.65727$ & $1.34273$\\
\parbox[0pt][2em][c]{0cm}{} & $\eta_{V[X],V[\cdot]}$ & $2.68546$ & $1.31454$\\
\parbox[0pt][2em][c]{0cm}{} & $\eta_{V,V[\cdot]}$ & $2.66666$ & $1.33333$\\
\hline
\end{tabular}
\end{center}

\paragraph{Edge-plate adjacencies.} Any edge is adjacent to exactly $3$ plates of the STIT tessellation, so $\la_P\mu_{P,E}=3\la_E$. Of course this holds for all the considered subclasses $\sE[\cdots]$, which leads to $$\mu_{P,E[\cdots]}={36\over 7}\eps_{E[\cdots]}.$$ With the $\eps$-values given by Theorems \ref{thm:TXclassification}-\ref{thm:Z1classification} we find the values summarized in the tables below.
\begin{center}
\begin{tabular}{|cc||c|c|c|}
\hline
\parbox[0pt][2em][c]{0cm}{} & & $E[TT]$ & $E[XX]$ & $E[TX]$\\
\hline
\parbox[0pt][2em][c]{0cm}{} & $\mu_{P,\cdot}$ & $2.27766$ & $0.56337$ & $2.30182$\\
\hline
\end{tabular}
\end{center}
\begin{center}
\begin{tabular}{|cc||c|c|c||c|c|c|}
\hline
\parbox[0pt][2em][c]{0cm}{} & & $E{[P_1,1]}$ & $E{[P_1,2]}$ & $E{[P_1,3]}$ & $E{[Z_1,0]}$ & $E{[Z_1,1]}$ & $E{[Z_1,2]}$\\
\hline
\parbox[0pt][2em][c]{0cm}{} & $\mu_{P,\cdot}$ & $2.85452$ & $1.54624$ & $0.74209$ & $3.21197$ & $1.18879$ & $0.74209$\\
\hline
\end{tabular}
\end{center}

\paragraph{Edge-facet adjacencies.} Any edge is adjacent to $5$ cell facets (elements from $\sZ_2$) and is located on the boundary of $4$ and in the relative interior of another one so that $\la_{Z_2}\mu_{Z_2,E}=5\la_E$, $\la_{Z_2}\mu_{\partial Z_2,E}=4\la_E$ and $\la_{Z_2}\mu_{\overset{\circ}{Z_2},E}=\la_E$, where the typical object of class $\sE$ can be replaced by the typical object of each of the subclasses of $\sE$ introduced in Subsection \ref{subsec:edgeclassification}. With $\lambda_E=2{\lambda_{Z_2}}$ this leads to
$$
 \mu_{Z_2,E[\cdots]} = 10\eps_{E[\cdots]},\qquad
\mu_{\overset{\circ}{Z_2},E[\cdots]} = 2\eps_{E[\cdots]},\qquad
 \mu_{\partial Z_2,E[\cdots]} = 8\eps_{E[\cdots]}.$$
The related mean values are with the $\eps$-values from Theorem \ref{thm:TXclassification}-\ref{thm:Z1classification} summarized below.
\begin{center}
\begin{tabular}{|cc||c|c|c|}
\hline
\parbox[0pt][2em][c]{0cm}{} & & $E[TT]$ & $E[XX]$ & $E[TX]$\\
\hline
\parbox[0pt][2em][c]{0cm}{} & $\mu_{Z_2,\cdot}$ & $4.42878$ & $1.09545$ & $4.47577$\\
\parbox[0pt][2em][c]{0cm}{} & $\mu_{\overset{\circ}{Z_2},\cdot}$ & $0.88576$ & $0.21909$ & $0.89515$\\
\parbox[0pt][2em][c]{0cm}{} & $\mu_{\partial Z_2,\cdot}$ & $3.54303$ & $0.87636$ & $3.58062$\\
\hline
\end{tabular}
\end{center}
\begin{center}
\begin{tabular}{|cc||c|c|c||c|c|c|}
\hline
\parbox[0pt][2em][c]{0cm}{} & & $E[P_1,1]$ & $E[P_1,2]$ & $E[P_1,3]$ & $E[Z_1,0]$ & $E[Z_1,1]$ & $E[Z_1,2]$\\
\hline
\parbox[0pt][2em][c]{0cm}{} & $\mu_{Z_2,\cdot}$ & $5.55046$ & $3.00657$ & $1.44296$ & $6.24550$ & $2.31154$ & $1.44296$\\
\parbox[0pt][2em][c]{0cm}{} & $\mu_{\overset{\circ}{Z_2},\cdot}$ & $1.11009$ & $0.60132$ & $0.28859$ & $1.24910$ & $0.46231$ & $0.28859$\\
\parbox[0pt][2em][c]{0cm}{} & $\mu_{\partial Z_2,\cdot}$ & $4.44037$ & $2.40526$ & $1.15437$ & $4.99640$ & $1.84923$ & $1.15437$\\
\hline
\end{tabular}
\end{center}

\paragraph{Edge-I-polygon adjacencies.} We observe that any edge belongs to exactly $2$ I-polygons. In one it is located on the boundary and in the relative interior of the other one so that $\la_I\mu_{I,E}=2\la_E$, $\la_I\mu_{\overset{\;\circ}{I},E}=\la_E$ and $\la_I\mu_{\partial I,E}=\la_E$, where the same equations also hold true for any of the edge classes from Subsection \ref{subsec:edgeclassification}. Hence, we obtain with $\lambda_E=12{\lambda_{I}}$ 
$$\mu_{I,E[\cdots]}=24\eps_{E[\cdots]},\qquad \mu_{\overset{\;\circ}{I},E[\cdots]}=12\eps_{E[\cdots]} \qquad \mu_{\partial I,E[\cdots]}=12\eps_{E[\cdots]}$$ and the following numerical values.
\begin{center}
\begin{tabular}{|cc||c|c|c|}
\hline
\parbox[0pt][2em][c]{0cm}{} & & $E[TT]$ & $E[XX]$ & $E[TX]$\\
\hline
\parbox[0pt][2em][c]{0cm}{} & $\mu_{I,\cdot}$ & $10.6291$ & $2.6291$ & $10.7418$\\
\parbox[0pt][2em][c]{0cm}{} & $\mu_{\overset{\;\circ}{I}/\partial I,\cdot}$ & $5.3145$ & $1.3145$ & $5.3709$\\
\hline
\end{tabular}
\end{center}
\begin{center}
\begin{tabular}{|cc||c|c|c||c|c|c|}
\hline
\parbox[0pt][2em][c]{0cm}{} & & $E[P_1,1]$ & $E[P_1,2]$ & $E[P_1,3]$ & $E[Z_1,0]$ & $E[Z_1,1]$ & $E[Z_1,2]$\\
\hline
\parbox[0pt][2em][c]{0cm}{} & $\mu_{I,\cdot}$ & $13.3211$ & $7.2158$ & $3.4631$ & $14.9892$ & $5.5477$ & $3.4631$\\
\parbox[0pt][2em][c]{0cm}{} & $\mu_{\overset{\;\circ}{I}/\partial I,\cdot}$ & $6.6606$ & $3.6079$ & $1.7316$ & $7.4946$ & $2.7738$ & $1.7316$\\
\hline
\end{tabular}
\end{center}

\paragraph{Edge-cell adjacencies.} Similar as above we observe that any tessellation edge is adjacent to exactly $3$ cells so that $\la_Z\mu_{Z,E}=3\la_E$. Moreover, we notice that any edge is located in the edge skeleton $sk(Z)$ of exactly $2$ cells, which implies $\la_Z=\mu_{sk(Z),E}=2\la_E$. The same equations are also true if $E$ is replaced by any considered subclass $E[\cdots]$. This yields the following relations: 
$$\mu_{Z,E[\cdots]}=36\eps_{E[\cdots]}, \qquad \mu_{sk(Z),E[\cdots]}=24\eps_{E[\cdots]}$$ 
and the following mean values.
\begin{center}
\begin{tabular}{|cc||c|c|c|}
\hline
\parbox[0pt][2em][c]{0cm}{} & & $E[TT]$ & $E[XX]$ & $E[TX]$\\
\hline
\parbox[0pt][2em][c]{0cm}{} & $\mu_{Z,\cdot}$ & $15.9436$ & $3.9436$ & $16.1128$\\
\parbox[0pt][2em][c]{0cm}{} & $\mu_{sk(Z),\cdot}$ & $10.6291$ & $2.6291$ & $10.7418$\\
\hline
\end{tabular}
\end{center}
\begin{center}
\begin{tabular}{|cc||c|c|c||c|c|c|}
\hline
\parbox[0pt][2em][c]{0cm}{} & & $E[P_1,1] $ & $E[P_1,2]$ & $E[P_1,3]$ & $E[Z_1,0]$ & $E[Z_1,1]$ & $E[Z_1,2]$\\
\hline
\parbox[0pt][2em][c]{0cm}{} & $\mu_{Z,\cdot}$ & $19.9817$ & $10.8237$ & $5.1947$ & $22.4838$ & $8.3215$ & $5.1947$\\
\parbox[0pt][2em][c]{0cm}{} & $\mu_{sk(Z),\cdot}$ & $13.3211$ & $7.2158$ & $3.4631$ & $14.9892$ & $5.5477$ & $3.4631$\\
\hline
\end{tabular}
\end{center}

\section{Summary and discussion}\label{sec:summary}

The average combinatorial structure of spatial random STIT tessellations was in the focus of this paper. We have introduced several interesting classes of geometric objects that are determined by the tessellation and have recalled some basic mean values from earlier works. The main contribution of the present work is a detailed refined combinatorial analysis that involves various edge classifications. For example, we were able to explore the geometry of the neighbors of the typical vertex as well as the typical T- and the typical X-vertex. The related mean values involve the probabilities that the typical edge belongs to one of our introduced edge classes. These probabilities in turn arise from a detailed distributional analysis of the typical I-segment, which is also further developed in this paper.

The reader might wonder why our refined combinatorial structures are based only on edges and not on a similar classification of plates, which are the $2$-dimensional primitive objects. Clearly, a classification through the type of vertices is not very promising, because a plate can have arbitrary many corners or vertices on its boundary. However, a classification according to the equality relationship for  plates with facets $Z_2$ would be possible and would lead to a breakdown of $\sP$ in to the classes $\sP[Z_2,i]$ and to the related probabilities $\eps_{P[Z_2,i]}=\bP(P\in\sP[Z_2,i])$ for $i=0,1,2$. To calculate these probabilities, one would need to investigate the geometry of I-polygons (instead of I-segments). From \cite{ST} is is known that the distribution of the typical I-polygon can be written as a mixture of typical cell distributions of certain Poisson line tessellations. In order to calculate $\eps_{P[Z_2,i]}$ one would have to calculate the area distribution of the typical I-polygon, which -- in view of the last sentence -- reduces to the determination of the area distribution of a typical cell in a Poisson line tessellation. This, however, is a long standing open problem in stochastic geometry. This way a calculation of the probabilities $\eps_{P[Z_2,i]}$ is currently out of reach. For this reason we had to restrict our work to the case of vertices and edges, which are analytically tractable, because of the line intersection property of STIT tessellations, which tells us that the intersection of a STIT tessellation with a line induces a homogeneous Poisson point process on that line. Point processes of this type are very well known and we can make use of their properties to gain detailed information about I-segments and edges.

Let us finally remark that for some of the quantities considered in this paper one can derive the whole distribution and not only the mean value. However, these distributions appear to be complicated in general and therefore we have here restricted our attention to first-order properties. To give a simple example for a distribution let us consider the number $\nu_{E[XX]}$ of $E[XX]$-edges that are located on the typical I-segment. It holds that $${\Bbb P}(\nu_{E[XX]}=0)=p_0+p_1+\sum_{n=2}^\infty(1-p_{X|n}^2)^{n-1}p_n$$ and that $${\Bbb P}(\nu_{E[XX]}=N)=\sum_{n=N+1}^\infty{n-1\choose N}p_{X|n}^{2N}(1-p_{X|n}^2)^{n-1-N},\qquad N\geq 1,$$ with $p_n$ and $p_{X|n}=1-p_{T|n}$ given by (\ref{eq:pm}) and (\ref{eq:ptm}), respectively. For small $N$ these expressions can be evaluated numerically and we find
\begin{alignat}{4}
\nonumber & {\Bbb P}(\nu_{E[XX]}=0) &&\approx 0.795161,\qquad && {\Bbb P}(\nu_{E[XX]}=1) &&\approx 0.131115,\\
\nonumber & {\Bbb P}(\nu_{E[XX]}=2) &&\approx 0.046467,\qquad && {\Bbb P}(\nu_{E[XX]}=3) &&\approx 0.016359.
\end{alignat}
For the other edge types these probabilities are considerably more involved and for this reason omitted.

\section{Proofs}\label{sec:proofs}

\subsection{More distributions for the typical I-segment}\label{subsec:subsec:isegmentdistributions2}

In this subsection we continue the distributional analysis of the typical I-segment started in Subsection \ref{subsec:isegmentdistributions1} and recall some other results from \cite{TWN}. They prepare the proofs presented below. Firstly, formula (\ref{eq:pm}) can be refined by considering the joint distribution of the number of T- and X-vertices in the relative interior of the typical I-segment. Let us define $p_{m,j}$ ($m,j\in\NN$) to be the probability that in its relative interior the typical I-segment has exactly $m$ T-vertices and $j$ X-vertices, i.e. $p_{m,j}=\bP(\nu_T=m,\nu_X=j)$. Then
\begin{equation}\label{eq:pXT}
p_{m,j}=3\cdot 2^m{m+j\choose m}\int\limits_0^1\int\limits_0^1 (1-a)^3a^m{(1-(1-a)(1-b))^j\over(3-(1-a)(2-b))^{m+j+1}}dbda,
\end{equation}
see \cite{TWN}.
Besides the two probabilities $p_n$ from (\ref{eq:pm}) and $p_{m,j}$, we also need the probability $p_{l,r}^{LR}$ ($l,r\in\NN$) that on the typical I-segment exactly $l$ internal T-vertices are induced from polygons pointing to the left and $r$ to the right. It is given by
\begin{equation}\label{eq:pLR}
p_{l,r}^{LR}=3{l+r\choose l}\int\limits_0^1(1-a)^3{a^{l+r}\over(1+a)^{l+r+1}}da.
\end{equation}
\paragraph{Proof of (\ref{eq:pLR}).} Let us consider the spatio-temporal construction of STIT tessellations in the isotropic case, that is the case where the plane measure $\L$ is the isometry invariant measure on $\cH$, see Subsection \ref{subsec:STITs}. (This is only for simplicity, the anisotropic case can also be considered, but then also the direction of the segment has to be taken into account. However, it appears that (\ref{eq:pLR}) is independent of the choice of the direction and hence independent of the plane measure $\L$ used in the STIT construction, see Subsection \ref{subsec:STITs}.) In this construction every I-segment receives a birth-time, which is defined as the birth-time of the I-polygon that creates the segment by intersection of another I-polygon that has been born earlier. The birth-time of the typical I-segment is denoted by $\beta$. Moreover, we denote by $\ell$ the length of the typical I-segment. From \cite{TWN} we know that the joint density $f_{\ell,\beta}(x,s)$ of $(\ell,\beta)$ for the STIT tessellation $Y(t)$ equals 
$$f_{\ell,\beta}(x,s)={3s^3\over 2t^3}e^{-{1\over 2}sx}\qquad x>0,\; 0<s<t.$$ We condition now on the event that $(\ell,\beta)=(x,s)\in(0,\infty)\times(0,t)$ and conclude from the intersection property of STIT tessellations that under this condition the number $\nu_L$ of left-pointing T-vertices and the number $\nu_R$ of right-pointing T-vertices in the relative interior of the typical I-segment are independent and Poisson distributed with parameter ${1\over 2}x(t-s)$, which is to say $$\bP(\nu_L=l,\nu_R=r|(\ell,\beta)=(x,s))={\left({1\over 2}x(t-s)\right)^{l+r}\over l!r!}e^{-x(t-s)},$$ see for example \cite{TWN} or the references cited therein. Averaging with respect to the joint density $f_{\ell,\beta}(x,s)$ yields
\begin{eqnarray}
\nonumber p_{l,r}^{LR} &=& \int\limits_0^t\int\limits_0^\infty \bP(\nu_L=l,\nu_R=r|(\ell,\beta)=(x,s))f_{\ell,\beta}(x,s)dxds\\
\nonumber &=& \int\limits_0^t\int\limits_0^\infty {\left({1\over 2}x(t-s)\right)^{l+r}\over l!r!}e^{-x(t-s)}\cdot {3s^3\over 2t^3}e^{-{1\over 2}sx}dxds\\
\nonumber &=& \int\limits_0^t{s^2\over t^3}{\left({1\over 2}(t-s)\right)^{l+r}\over l!r!}{(l+r)!\over\left({1\over 2}(2t-s)\right)^{l+r+1}}ds\\
\nonumber &=& 3{l+r\choose l}\int\limits_0^t{s^2\over t^3}{(t-s)^{l+r}\over(2t-s)^{l+r+1}}ds\\
\nonumber &=& 3{l+r\choose l}\int\limits_0^1(1-a)^3{a^{l+r}\over(1+a)^{l+r+1}}da,
\end{eqnarray}
where in the last step we have applied the substitution $s=t(1-a)$. This proves our claim.\hfill $\Box$

\subsection{Proofs for the results of Subsection \ref{subsec:isegmentdistributions1}}\label{subsec:proofth1}

\paragraph{Proof of Theorem \ref{thm:vertexinisegment} (a).} Recall that by $p_{m,j}$ ($m,l\in\NN$) we denote the probability that the typical I-segment contains exactly $m$ internal vertices from $\sV[T]$ and $j$ internal vertices from $\sV[X]$. Note, moreover, that $1-p_0=1-\bP(\nu=0)$ is the probability that the typical I-segment has internal vertices. Thus, given the typical I-segment contains exactly $n\geq 1$ internal vertices, and given that $m\in\{0,\ldots,n\}$ of them are of type $T$, the probability for choosing a T-vertex is $m/n$. Averaging over all possible numbers $m$ and $n$ we find 
\begin{equation}\label{EqpT}
p_T={1\over 1-p_0}\sum_{n=1}^\infty\sum_{m=0}^n{m\over n}p_{m,n-m},
\end{equation}
where the factor $1/(1-p_0)$ comes from conditioning on the event that the typical I-segment has internal vertices as explained above. Combining now (\ref{EqpT}) with (\ref{eq:pXT}) we calculate $p_T$ as follows:
\begin{eqnarray}
\nonumber p_T &=& {1\over 1-p_0}\sum_{n=1}^\infty\sum_{m=0}^n{m\over n}p_{m,n-m}\\
\nonumber &=& {3\over 1-p_0}\int\limits_0^1\int\limits_0^1\sum_{n=1}^\infty\sum_{m=0}^n{m\over n}{n\choose m}\\
\nonumber & & \hspace{3.3cm}\times(2a)^m (1-a)^3{(1-(1-a)(1-b))^{n-m}\over(3-(1-a)(2-b))^{n+1}}dbda\\
\nonumber &=& {3\over 1-p_0}\int\limits_0^1\int\limits_0^1{2a(1-a)^2\over 1+2a+b-ab}dbda={3\over 1-p_0}\left({7\over 2}+{28\over 3}\ln 2-9\ln 3\right).
\end{eqnarray}
The precise value for $p_T$ is obtained by taking into account that $p_0={189\over 8}\ln 3-26\ln 2-{15\over 2}$. In addition, the relation $p_T+p_X=1$ implies the value for $p_X$, which completes our argument.\hfill $\Box$

\paragraph{Proof of Theorem \ref{thm:vertexinisegment} (b).} Now we use for the internal vertices of type T on the typical I-segment the same method as in the proof of part (a) above and show $p_{L|T}=p_{R|T}=\frac{1}{2}$, see Remark 2 following Theorem \ref{thm:vertexinisegment}. Given the typical I-segment contains m internal vertices of type T, and $l$ of them point to the left, the probability for choosing a left-pointing T-vertex is $l/m$. Averaging yields for $p_{L|T}$ the expression $$p_{L|T}={1\over \bP(\nu_T\geq 1)}\sum_{m=1}^\infty\sum_{l=0}^m{l\over m}p_{l,m-l}^{LR},$$ 
where the factor 1/$\bP(\nu_T\geq 1)$ comes from conditioning on the event $\nu_T\geq 1$. To evaluate the double sum we use (\ref{eq:pLR}), which yields ${17\over 2}-12\ln 2$ by interchanging summation with integration, compare with the proof of part (a). It remains to calculate the probability $\bP(\nu_T\geq 1)$, which by using (\ref{eq:pXT}), is given by $$\bP(\nu_T\geq 1)=\sum_{m=1}^\infty\sum_{j=0}^\infty p_{m,j}=17-24\ln 2.$$ We thus find $$p_{L|T}={{17\over 2}-12\ln 2\over 17-24\ln 2}={1\over 2}$$ and $p_{R|T}=1-p_{L|T}={1\over 2}$.\hfill $\Box$

\subsection{Proofs for the results of Subsection \ref{subsec:edgeclassification}}

\paragraph{Proof of Theorem \ref{thm:TXclassification}.} We start by noting that any edge is located on exactly one I-segment and that the two endpoints of an I-segment are T-vertices. Now we consider on the typical I-segment the different edge-types induced by the different internal vertices. If $\nu=0$ (the typical I-segment has no internal vertices) we have exactly one edge of type TT. If $\nu=1$, then the typical I-segment comprises two TT-edges, if the internal vertex is of type T -- an event having probability $p_{T|1}$. And we have two TX-edges, if it is an X-vertex (with probability $p_{X|1}$). For $\nu \geq 2$ at first we are interesting in the edges, whose endpoints both are internal vertices of the I-segment. For $\nu=n$ there are $n-1$ of these edges on the typical I-segment. If we uniformly select one of them, it is a TT-edge with probability $p^2_{T|n}$, thanks to the independence of the vertex type of two neighboring vertices on an I-segment. And in the same way it is an XX-edge with probability $p^2_{X|n}$ and a TX-edge with probability $2p_{T|n}p_{X|n}$. Furthermore, there are two ``boundary'' edges (where one endpoint is also an endpoint of the I-segment and therefore a T-vertex), they are of type TT with probability $p_{T|n}$, of type TX with probability $p_{X|n}$ and of type XX with probability $0$.

For this reason the mean number $\mu_{{I_1},E[TT]}$ of edges of type TT that are located on the typical I-segment of the tessellation can be calculated by 
\begin{equation}\label{proofthm1eq1}
\mu_{{I_1},E[TT]}=p_0+2p_{T|1}p_1+\sum_{n=2}^\infty(2p_{T|n}+(n-1)p_{T|n}^2)p_n. 
\end{equation}
We observe next that the intensity $\la_{E[TT]}$ of TT-edges equals $\la_{E[TT]}=\la_{I_1}\mu_{{I_1},E[TT]}$. Moreover, we have the relation $\la_E=3\la_{I_1}$ and thus the probability $\eps_{E[TT]}$ is given by $$\eps_{E[TT]}={\la_{E[TT]}\over\la_E}={1\over 3}\mu_{{I_1},E[TT]},$$ which in view of (\ref{proofthm1eq1}) proves our first claim. 

To calculate the intensities $\la_{E[XX]}$ and $\la_{E[TX]}$ of XX- and TX-edges we observe that any vertex (regardless of its type) has exactly $4$ outgoing edges, which implies the two intensity relationships $$2\la_{E[XX]}+\la_{E[TX]}=4\la_{V[X]}\quad{\rm and}\quad 2\la_{E[TT]}+\la_{E[TX]}=4\la_{V[T]}$$ by counting the X- and T-vertices through the different types of edges. Thus, $$\eps_{E[TX]}={\la_{E[TX]}\over\la_E}={4\la_{V[T]}-2\la_{E[TT]}\over\la_E}=4{\la_{V[T]}\over\la_E}-2\eps_{\sE[TT]}={4\over 3}-2\eps_{E[TT]}$$ and similarly 
$\eps_{E[XX]}=\eps_{E[TT]}-{1\over 3}$, which completes the proof.\hfill $\Box$

\begin{remark}
Alternatively, the probabilities $\eps_{E[XX]}$ and $\eps_{E[TX]}$ are determined by $$\eps_{E[XX]}={1\over 3}\mu_{{I_1},E[XX]}\qquad{\rm and}\qquad\eps_{E[TX]}={1\over 3}\mu_{{I_1},E[TX]},$$ where $\mu_{{I_1},E[XX]}$ and $\mu_{{I_1},E[TX]}$ are, respectively, the mean number of edges of type XX and TX located on the typical I-segment. Using the same argument as in the proof above, one can show that $\mu_{{I_1},E[XX]}$ and $\mu_{{I_1},E[TX]}$ satisfy $$\mu_{{I_1},E[XX]}=\sum_{n=2}^\infty (n-1)p_X^2p_n$$ and $$\mu_{{I_1},E[TX]}=2p_{X|1}p_1+\sum_{n=2}^\infty(2p_{X|n}+2(n-1)p_{X|n}p_{T|n})p_n,$$ which leads to the same values.
\end{remark}

\paragraph{Proof of Theorem \ref{thm:P1classification}.} This can be shown with the help of the same technique as already used in the proof of Theorem \ref{thm:TXclassification} above. We first find for the mean number of $E[P_1,i]$-edges that are located in the relative interior of the typical I-segment
$$\eps_{E[P_1,i]}={\la_{E[P_1,i]}\over\la_E}={1\over 3}\mu_{{I_1},E[P_1,i]}, \ \  i=1,2,3 .$$

However, an edge belongs to class $\sE[P_1,3]$ if and only if it is an I-segment at the same time. Hence, $\mu_{{I_1},E[P_1,3]}=\bP(\nu=0)=p_0$, which can be evaluated using (\ref{eq:pm}). In addition, if $\nu=1$ both edges on the I-segment are of type $E[P_1,2]$. 

For $\nu\geq 2$ the two ``boundary'' edges are again in $\sE[P_1,2]$. An edge whose both endpoints are internal vertices of the I-segment is of type $E[P_1,2]$ if its endpoints are X and X, L and L or R and R, see Figure \ref{fig:edges}. It is of type $E[P_1,1]$ if its endpoint are of type T and X or L and R, respectively.
So, we find $$\mu_{{I_1},E[P_1,2]}=2p_1+\sum_{n=2}^\infty\left(2+(n-1)p_{T|n}^2/2+(n-1)p_{X|n}^2\right)p_n$$ and similarly $$\mu_{{I_1},E[P_1,1]}=\sum_{n=2}^\infty\left((n-1)p_{T|n}^2+2(n-1)p_{T|n}p_{X|n}\right)p_n,$$ which completes the argument. \hfill $\Box$ 

\paragraph{Proof of Theorem \ref{thm:Z1classification}.} Because this follows once again by similar arguments as above we restrict ourself to a rough sketch. We have again $$\eps_{E[Z_1,j]}={\la_{E[Z_1,j]}\over\la_E}={1\over 3}\mu_{{I_1},E[Z_1,j]},\qquad j=0,1,2,$$ where $\mu_{{I_1},E[Z_1,j]}$ is the mean number of $E[Z_1,j]$-edges in the relative interior of the typical I-segment and hence it remains to determine $\mu_{{I_1},E[Z_1,i]}$. 

However, as in the proof of Theorem \ref{thm:P1classification}, $\mu_{{I_1},E[Z_1,2]}=\bP(\nu=0)=p_0$ and further $$\mu_{{I_1},E[Z_1,1]}=2p_{T|n}p_1+\sum_{n=2}^\infty\left(2p_{T|n}+(n-1)p_{T|n}^2/2\right)p_n$$ and
\begin{eqnarray}
\nonumber \mu_{{I_1},E[Z_1,0]} &=& 2p_{X|n}p_n+\sum_{n=2}^\infty\left(2p_{X|n}+(n-1)p_{X|n}^2+2(n-1)p_{T|n}p_{X|n}\right.\\
\nonumber & & \hspace{3.5cm}\left. +(n-1)p_{T|n}^2/2\right)p_n,
\end{eqnarray}
which may be seen by taking into account the different possibilities for subsequent internal vertices in the relative interior of the typical I-segment and the two ``boundary'' edges, and see again Figure \ref{fig:edges}.\hfill $\Box$

\subsection{Proof of the vertex-edge adjacencies in Subsection \ref{subsec:newmeanvalues}}\label{subsec:proofvertexedge}
To establish the mean value relations for the mean number of edges equal to one $P_1$-segment adjacent to the typical T- or X-vertex $\mu_{V[\cdot],E[P_1,1]}$ first we consider the edges in subclass $\sE[P_1,1]$. All of them are edges on an I-segment whose both endpoints are internal vertices of the I-segment. Moreover, these endpoints must be of type L and R or of type T (L or R) and X, see Figure \ref{fig:edges}. Using the same method as in the previous proofs we obtain
$$\lambda_{V[T]}\mu_{V[T],E[P_1,1]} = \lambda_{I_1} \sum_{n=2}^{\infty} (n-1) (p^2_{T|n} + p_{T|n}p_{X|n})p_n.$$
With $\lambda_{I_1}=\frac{2}{3}{\lambda_V}$, $\lambda_{V[T]}=\frac{2}{3}{\lambda_V}$ and $p_{X|n}=1-p_{T|n}$ we thus find
$$\mu_{V[T],E[P_1,1]}= \sum_{n=2}^{\infty} (n-1)  p_{T|n}p_n.$$ The sum can be evaluated explicitly by using (\ref{eq:pXT}), which yields
\begin{eqnarray}
\nonumber \mu_{V[T],E[P_1,1]} &=& \sum_{n=2}^\infty(n-1)\sum_{k=0}^n{k\over n}2^k{n\choose k}\\
\nonumber & & \hspace{1cm}\times\int\limits_0^1\int\limits_0^1a^k(1-a)^3{1-(1-a)(a-b)^{n-k}\over(3-(1-a)(2-b))^{n+1}}dbda\\
\nonumber &=& \int\limits_0^1\int\limits_0^1{6a(1-a)(a(3-b)+b)\over 1+a(2-b)+b}dbda\\
\nonumber &=& 27\ln 3-28\ln 2-{19\over 2}.
\end{eqnarray}
Similar considerations for the typical X-vertex imply
$$\lambda_{V[X]}\mu_{V[X],E[P_1,1]} = \lambda_{I_1} \sum_{n=2}^{\infty} (n-1)  p_{T|n}p_{X|n}p_n$$
and with ${\lambda_{I_1}}={2\over 3}\lambda_V$ and $\lambda_{V[X]}=\frac{1}{3}\lambda_V$ this reduces to
$$\mu_{V[X],E[P_1,1]}= 2 \sum_{n=2}^{\infty} (n-1)  p_{T|n}p_{X|n}p_n.$$
The relations for $\mu_{V[\cdot],E[P_1,2]}$ for both types of typical vertices follow from $\mu_{V[\cdot],E[P_1,1]}+\mu_{V[\cdot],E[P_1,2]}+\mu_{V[\cdot],E[P_1,3]}=4.$
\hfill $\Box$

\subsection*{Acknowledgements}
We feel grateful to Claudia Redenbach (Kaiserslautern) for providing the pictures of the STIT tessellations shown in Figure \ref{fig:STITconstruction} and Figure \ref{fig:STITconstructionfinal}.\\ The second author was supported by the German research foundation (DFG), grant WE 1799/3--1.



\end{document}